\documentclass[12pt]{amsart}

\usepackage{amssymb,amsmath,graphicx,color,amsthm, amsfonts, mathbbol}
 
\definecolor{Red}{cmyk}{0,1,1,0}

\definecolor{verde}{cmyk}{1,0,1,0}

\definecolor{azul}{cmyk}{1,1,0,0}

\numberwithin{equation}{section}

\newcommand{\dmax}[1]{d_{\text{max}}(T_{#1})}
\newcommand{\ignore}[1]{}

\def\Ed{{\mathbb{E}}}
\def\Pd{{\mathbb{P}}}

\newcommand{\N}{\mathbb{N}}

\newcommand{\R}{\mathbb{R}}
\newcommand{\Z}{\mathbb{Z}}

\newcommand{\Q}{\mathbb{Q}}

\newcommand{\e}{\varepsilon}

\newcommand{\s}{\sigma}

\newcommand{\D}{\Delta}

\newcommand{\Ind}[1]{\mathbb{1}_{#1}}
\def\sA{\mathcal{A}}

\newcommand{\dist}[3]{\mathrm{dist}_{T_{#1}}(#2,#3)}
\newcommand{\degr}[2]{d_{T_{#1}}(#2)}

\def\Prob{{\rm Prob}}

\newcommand{\be}{\begin{equation}}
\newcommand{\ee}{\end{equation}}

\newtheorem{theorem}{Theorem}[section]
\newtheorem{remark}{Remark}[section]
\newtheorem{example}{Example}[section]
\newtheorem{proposition}[equation]{Proposition}
\newtheorem{definition}{Definition}[section]
\newtheorem{lemma}{Lemma}[section]
\newtheorem{corollary}{Corollary}[section]

\newtheorem{problem}{Problem}[section]

\newtheorem{obs}{Observation}[section]
\newtheorem{claim}{Claim}[section]
\newtheorem{fact}{Fact}[section]

\newenvironment{claimproof}[1]{\par\noindent\textit{Proof of the claim:}\space#1}{\hfill $\blacksquare$}

\def\sG{\mathcal{G}}
\def\sT{\mathcal{T}}

\begin{document}
\title{Building your path to escape from home}
\author{D. Figueiredo$^1$}

\author{G. Iacobelli$^1$}

\author{R. Oliveira$^2$}

\author{B. Reed$^{3 \; 4}$}
\author{R. Ribeiro$^2$}

\date{\today \\
	$^1$ UFRJ, Universidade Federal do Rio de Janeiro\\
	$^2$ IMPA, Instituto de Matem\'atica Pura e Aplicada \\
	$^3$ School of Computer Science, McGill Univeristy \\
	$^4$ Laboratoire I3S, CNRS, France.
}

\begin{abstract}
Random walks on dynamic graphs have received increasingly more attention from different academic communities over the last decade. Despite the relatively large literature, little is known about random walks that construct the graph where they walk while moving around.
In this paper we study one of the simplest conceivable discrete time models of this kind, which works as follows: before every walker step, with probability $p$ a new leaf is added to the vertex currently occupied by the walker. The model grows trees and we call it the Bernoulli Growth Random Walk (BGRW). 
We show that the BGRW walker is transient and has a well-defined linear speed $c(p)>0$ for any $0<p\leq 1$. We also show that the tree as seen by the walker converges (in a suitable sense) to a random tree that is one-ended. Some natural open problems about this tree and variants of our model are collected at the end of the paper.
\end{abstract}
\maketitle

\section{Introduction}

Random walks on graphs that change over time have received much attention over the past decades. Within this context, a large body of work assumes only edge (or node) weights change over time while the graph structure (i.e., edge set) remains constant. Examples include reinforced random walks and random walks in random environments~\cite{Amir_Changing,cotar2017edge,Huang_EvolvingSets,disertori2015transience,kious2017phase}. A much smaller line of work assumes that the graph structure (edges and nodes) changes over time. However, such works generally assume graph dynamics to be independent of the walker~\cite{avena2016mixing,figueiredo2012,peres2015random} (an exception is~\cite{iacobelli2016}). 

This work explores a novel model where the random walk constructs its own graph, mutually coupling the walker and graph dynamics. This model is defined as follows: 
\begin{enumerate}
	\item[(0)] start with a finite tree with the walker sitting on one of the vertices; 
	\item[(1)] with probability $p$, add and connect a new leaf to the current location of the walker; 
	\item[(2)] let the walker take one step on the current graph;
	\item[(3)] go to step 1.
\end{enumerate}
We refer to this model as BGRW (Bernoulli Growth Random Walk). Note that or model may be seen as a sequence of pairs $(T_n,X_n)_{n \in \mathbb{N}}$, where $T_n$ is a tree and $X_n$ is the walker's position on the tree $T_n$. A more formal definition is given in Section~\ref{sec:model}. This model is a variant of the Non-Restart Random Walk  (NRRW) proposed by Amorim et al. ~\cite{Amorim_NRRW}. There, the initial graph is a vertex with a loop edge, and new leafs are created every $s$ steps for a fixed $s>0$. Prior to \cite{Amorim_NRRW}, other models of walks creating graphs had been proposed, but they all had periodic {\em restart} step where the random walker would jump to a uniform vertex in the tree~\cite{cannings2013,saramaki2004scale}. 

\begin{figure}[h]
\begin{center}
\scalebox{0.68}{
\begin{tabular}{|c|c|c|c|}
 \hline
\includegraphics[scale=0.27]{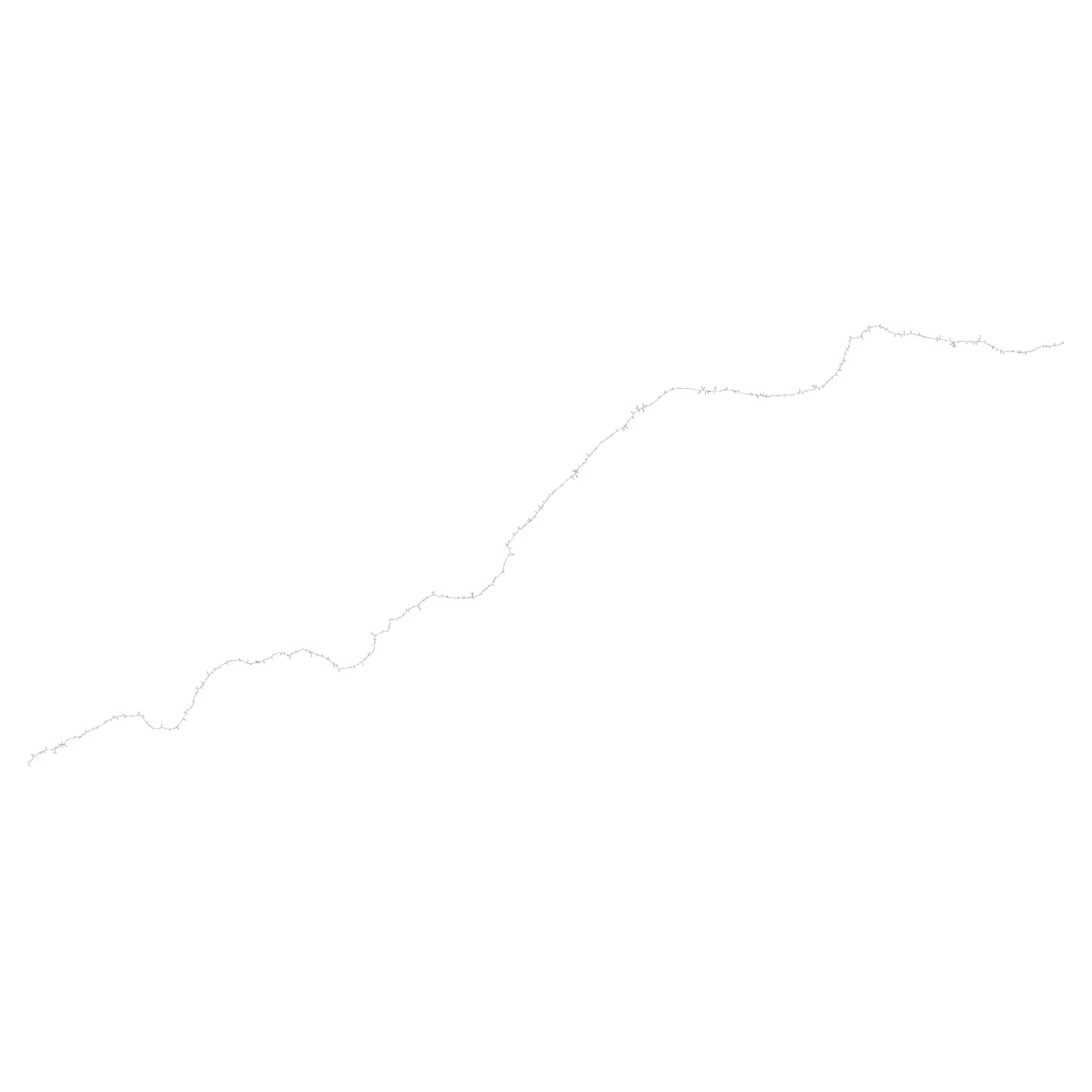}
 &
\includegraphics[scale=0.27]{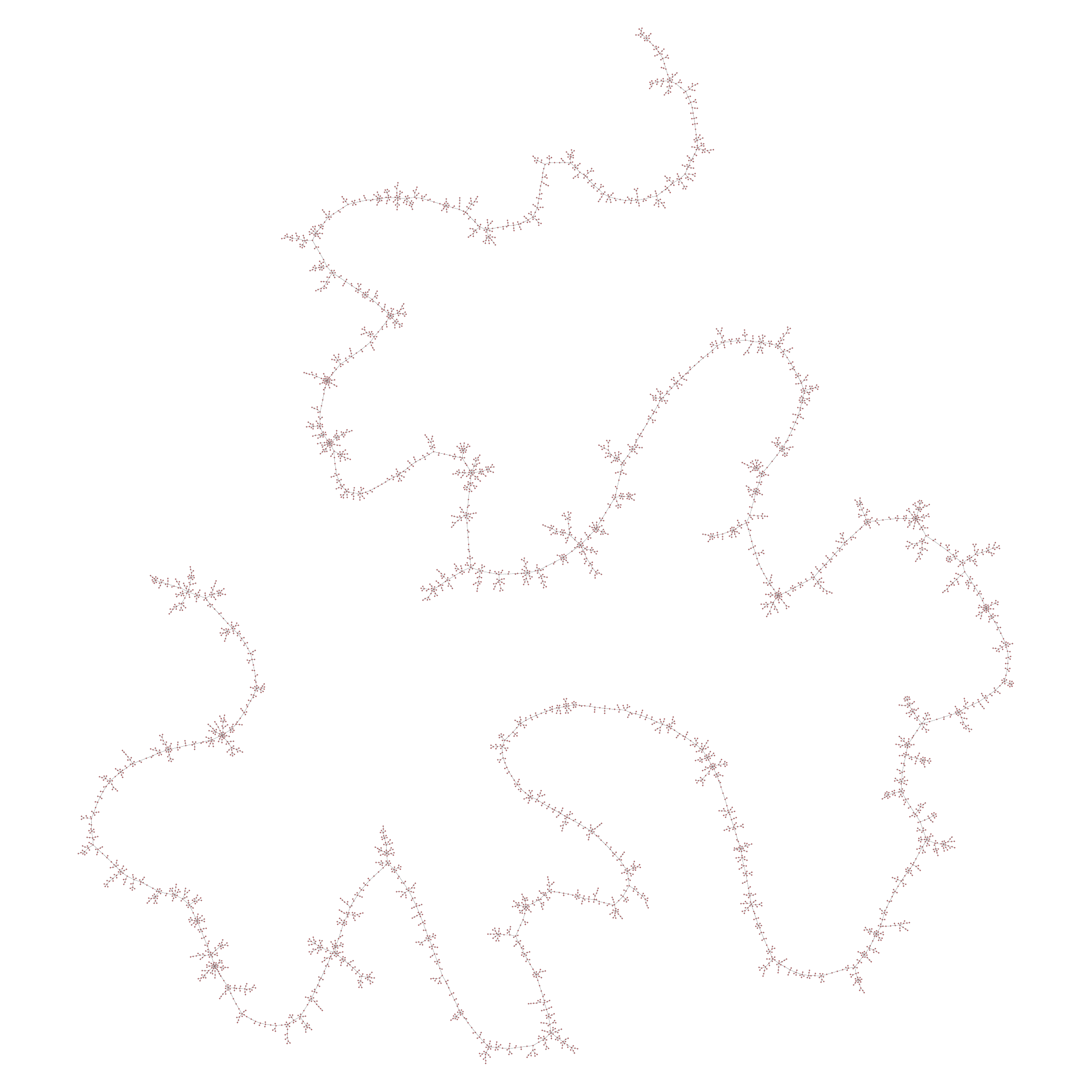}
 &
\includegraphics[scale=0.27]{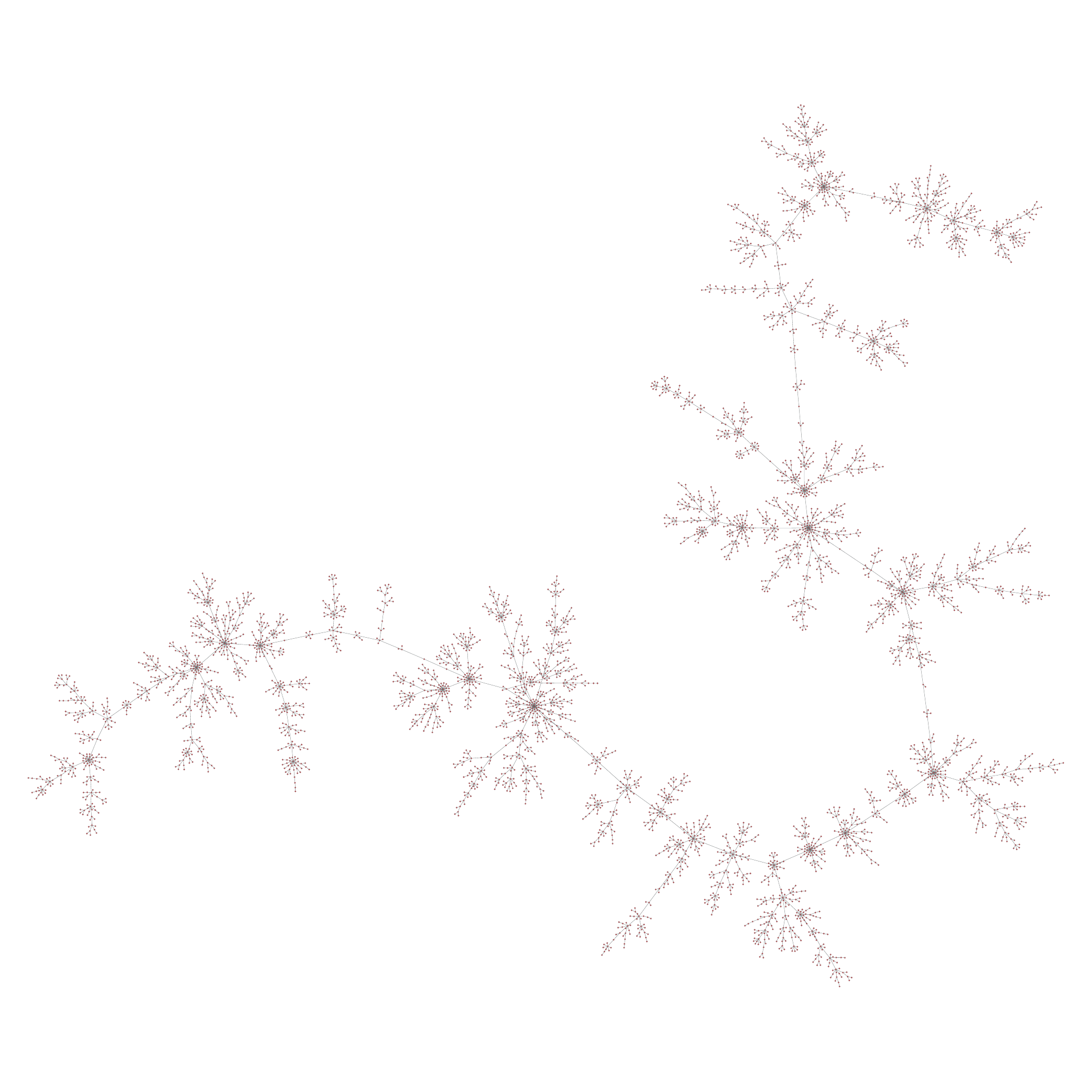}
&
\includegraphics[scale=0.27]{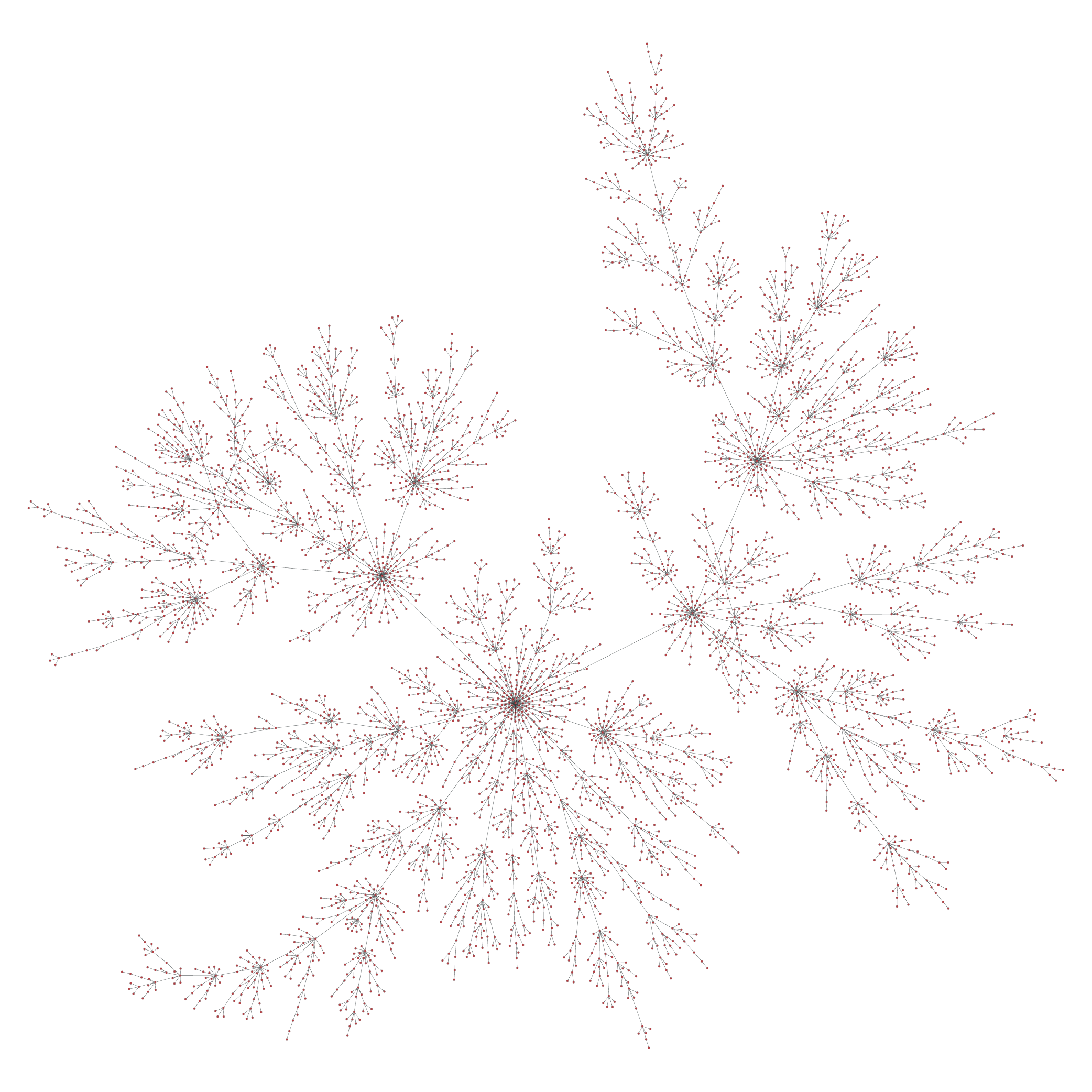}
 \\
\hline
\multicolumn{1}{c}{ } & \multicolumn{1}{c}{} & \multicolumn{1}{c}{} & \multicolumn{1}{c}{} \\
\multicolumn{1}{c}{\large $p=0.9$}
&
\multicolumn{1}{c}{\large  $p=0.5$}
&
\multicolumn{1}{c}{ \large $p=0.1$}
&
\multicolumn{1}{c}{ \large  $p=0.01$}
\end{tabular}
}
\caption{Trees with $n=5000$ nodes generated by simulating BGRW using different values for $p$, illustrating the diversity in the tree structure.}\label{fig:simul}
\label{fig:network}
\end{center}
\end{figure}

A fundamental question on dynamic graphs is the recurrent or transient nature of the walker. Figure~\ref{fig:simul} shows simulated sample paths from the BGRW model for different values of $p$. While all trees have exactly the same number of nodes ($n=5000$) their structural difference is striking. For larger values of $p$ the generated trees are very slim and long, as $p$ decreases the trees become fatter and shorter. Intuitively, with large $p$ the random walk can escape more easily, while for small $p$ the random walk wanders more. 

%

Our findings show that for any fixed $p$, the random walk escapes with a positive speed $c(p)>0$. This is similar to the results of Amorim et al. \cite{Amorim_NRRW} on the NRRW model with $s=1$; a similar result is expected to hold for all odd $s>1$ \cite{FIN17} \footnote{The case of even $s>1$ is recurrent due to parity issues \cite{FIN17}.}. Indeed, when $p$ is close to $1$, the proof for $s=1$ presented in \cite{Amorim_NRRW} can be adapted to show transience. 

As mentioned, most prior models and results constrain the walker to some graph. Our finding suggests that when the walker is {\em a prori} unconstrained, it constructs a graph that allows it to escape. In fact, we will see that the initial graph is in some sense ``forgotten". More specifically, we will see that the tree around the walker converges (in a suitable sense) to an infinite random tree whose law does not depend on the initial states.  We believe these findings are a significant contribution to the ongoing discussion of random walks on dynamic graphs.

\subsection{Main results}

Our first main result shows that our process has a well-defined, linear\footnote{We also use the term ``positive speed" to mean linear speed.} asymptotic speed. 

\begin{theorem}[Positive speed]\label{thm:positivespeed} For each $0<p\leq 1$, there exists a well-defined speed $0<c(p)\leq p$ for the BGRW process with parameter $p$. That is, for any initial condition~$(T_0,x_0)$ of the process, the distance between $X_n$ and $x_0$ satisfies:

\[\lim_{n \to \infty}\,\frac{{\rm dist}_{T_n}(X_n,x_0)}{n}= c(p), \;\mathbb{P}_{T_0,x_0,p}\mbox{-almost surely.}\]\end{theorem}

To obtain this result, we had to understand the distribution of degrees that $X_t$ sees along the way. A more general question is what is the frequency of (finite) trees that $X_t$ sees in a radius $r\geq 0$ around itself. To make this precise, we consider pairs $(T_t,X_t)$ as random elements of the space $\sT_*$ of {\em rooted locally finite trees} considered up to isomorphisms, with the {\em local topology}. The definition of $\sT_*$ and its local metric (which makes it a Polish space) are recalled in Section \ref{sec:finiteinfinite}. 

For each $t\in\N$, one may define an element $[T_t,X_t]\in\sT_*$ as the equivalence class of the tree ~$T_t$ rooted at $X_t$. In what follows, $[T_t,X_t]_r$ is the subtree of height $r$ consisting of $X_t$ and all nodes of $T_t$ within distance at most $r$ from $X_t$. The next theorem counts how many times the tree of radius $r$ around $X_t$ takes a certain shape. 

\begin{theorem}[Convergence of the tree as seen by the walker]\label{thm:localweak} For each $0<p\leq 1$, there exists a probability distribution $P_p$ over rooted trees such that, for any initial condition $(T_0,x_0)$ and any finite rooted tree $[S,x]\in\sT_*$ with height $r>0$,
\[\frac{1}{n+1}\sum_{t=0}^n\Ind{\{[T_t,X_t]_r=[S,x]\}}\rightarrow P_p(\{[T,o]\in\sT_*\,:\,[T,o]_r=[S,x]\})\;, \quad \mathbb{P}_{T_0,x_0,p}\mbox{-almost surely.}\]\end{theorem}

This is essentially a mixing property of the BGRW when viewed as a Markov chain over $\sT_*$. Our proof reveals that $P_p$ is the unique stationary measure of this chain with the following property: if the process is started from $P_p$, then almost surely, for all $\ell\geq 1$, the walker will eventually be at the tip of a path of length $\ell$. Theorem \ref{thm:localweak} follows from a stronger statement, Theorem \ref{thm:convergence} in Section \ref{sec:localweak}. 

The measure $P_p$ is an interesting object in itself. The next Theorem gives us some limited information on this probability measure. Recall that an infinite rooted tree is one-ended if it contains a single infinite path starting from the root. 

\begin{theorem}[Support of $P_p$]\label{thm:oneended} Given $0<p\leq 1$,  let $P_p$ be the limit measure in Theorem~\ref{thm:localweak}. Then $P_p$ is supported on one-ended infinite trees. Any rooted tree $[S,x]$ (with a certain height $r$) satisfies
\[P_p(\{[T,o]\in\sT_*\,:\,[T,o]_r=[S,x]\})>0.\]
Moreover, the degree of the root has exponential tails under $P_p$: that is, there exist $c,C>0$ such that for all $k\geq 1$:
\[P_p(\{[T,o]\in\sT_*\,:\,{\rm deg}_T(o)\geq k\})\leq C\,e^{-ck}.\]
\end{theorem}

Theorem \ref{thm:oneended} is in fact a corollary of Theorem \ref{thm:oneended2}, which we present in Section \ref{sec:localweak}.

\subsection{Intuition and some difficulties}

Behind our Theorems, there is a simple probabilistic mechanism that explains our claims. 

\begin{enumerate}
\item Given a constant $\ell$, the BGRW will often create new induced paths of length $\ell$. That is, there will be frequent times when the walker is at the ``tip"~of a path of length $\ell$. 
\item On the other hand, the probability of backtracking on a path of length $\ell$ by more than $\ell/2$ steps goes to zero very fast with $\ell$. 
\end{enumerate}

The moments when BGRW creates a long path and does not backtrack by more than half its length are ``local regeneration times" in the following precise sense: the sequence of neighborhoods of radius $\ell/2$ seen by the walker from that point on are independent of the past. From this perspective, the Laws of Large Numbers implicit in Theorems \ref{thm:positivespeed} and \ref{thm:localweak} are not surprising. 

As it turns out, there are difficulties in establishing the intuitive picture drawn by items $1$ and $2$ above. Creating paths looks easy enough: all one needs is to create a new leaf and jump to it $\ell$ consecutive times. However, in order for this to happen often, the walker needs to come across lots of nodes of low degree. For this to take place, it is necessary that nodes are usually not visited too many times. To prove that returns to a vertex are unlikely to be numerous, one needs to show that the walk moves away fast enough from any vertex. The upshot is that  a weaker form of {\em positive speed} is needed before we can justify the above items. 

Proving that ${\rm dist}_{T_n}(X_n,X_0)$ grows {\em at least linearly with time} will thus be the first major goal in our proof. Omer Angel (personal communication) suggested to us that this might be done comparing our process to a once reinforced random walk on trees. The rationale is that one may imagine that BGRW ``discovers" new edges rather than create them.  Unfortunately, we could not see a way to make the elegant arguments of Colevecchio \cite{collevecchio2006,collevecchio2017speed} work in our setting. We thus resort to a bare-hands argument.

Given the preliminary result on the growth of ${\rm dist}_{T_n}(X_n,X_0)$, we can move on to establishing the existence of the speed and the convergence of the tree as seen by the walker. It will be crucial for us that this tree -- or rather, the corresponding empirical measure -- converges almost surely to an invariant measure $P_p$ for the dynamics on rooted trees, for all ``nice"~initial conditions. There is a high-level similarity to the work of Lyons, Pemantle and Peres \cite{lyons1995} on random walk on Galton-Watson trees, where the existence of a speed relies on ergodic-theoretic arguments in the space $\sT_*$. However, in that paper the stationary measure can be described explicitly (leading to an explicit formula for the speed), which is not the case in our setting. Our analysis will also require more quantitative estimates on the convergence. Finally, we note that the random trees we consider are very different (eg. our tree is a.s. one-ended).

\subsection{Organization and main proof steps}

The remainder of the paper is organized as follows. In Section \ref{sec:model} we define our process formally. One important conceptual point will be to define it on locally finite trees from the start.

The proof of linear growth of $\dist{n}{X_n}{X_0}$ begins in Section \ref{sec:powersoflog}, where we show that, for any positive $M$, it is very likely that $\dist{m}{X_m}{X_0}\geq \log^M n$ for some $m\leq n$. This requires comparisons with simple one dimensional random walk. The argument continues in Section \ref{sec:loopprocess}, where we prove that backtracking on a long path, if it happens at all, typically takes a really long time. For this we introduce a simplified ``loop process" that only keeps track of the path itself along the way. Combining polylog distances and long times to backtrack, we prove in Section \ref{sec:positivedrift} that our process has positive drift away from $X_0$, and derive some consequences of this fact for the degrees in the tree. 

We switch gears for the remainder of the paper, as our arguments become more abstract. In Section \ref{sec:finiteinfinite} we extend the definition of our process to the space of rooted locally finite trees. The statement in Theorem \ref{thm:localweak}, on the convergence of the tree as seen by the walker, can then be stated as a weak convergence result for the empirical measure of this process. Tightness and weak convergence criteria are discussed in this section, and tightness is proven right away. 

Section \ref{sec:positivedrift} begins to study the loss-of-memory mechanism described above, whereby long paths are created and never backtracked on. Since we deal with infinite trees, we need some condition on the initial distribution to ensure that this takes place. We use this in Section~\ref{sec:finiteinfinite} to prove that the averages of ``local functions" along the trajectories of the process always converge almost surely. This is basically the last ingredient we need to prove stronger versions of Theorems \ref{thm:localweak} and \ref{thm:oneended} in Section \ref{sec:localweak}. One key point is that the limiting measure $P_p$ is a stationary distribution for our process. 

The paper wraps up with Section \ref{sec:final}, with some final comments, and an Appendix containing technical results.

\section{Definition of the model}\label{sec:model}
\subsection{Preliminaries} All trees in this paper are locally finite (all vertices have finite degree). Given a tree $T$, we let $V(T)$ and $E(T)$ denote its vertex and edge sets. For $x,y\in V(T)$,~$d_T(x)$ denotes the degree (number of neighbors) of $x$ in $T$ and ${\rm dist}_T(x,y)$ is the shortest-path distance between $x$ and $y$. 

We let $\Omega$ denote the set of all pairs $(T,x)$, where $T$ is locally finite tree with $V(T)\subset \N$ and~$\N\backslash V(T)$ infinite, and $x\in V(T)$. This set $\Omega$ can be described as a subset of a product space $\{0,1\}^{\N\cup \binom{\N}{2}}\times \N$ with a natural $\sigma$-field; we omit the details.  

Given a probability measure $\mu$ over some space, we use the symbol $U\sim \mu$ to mean that $U$ is a random element with law $\mu$.

\subsection{Definition of the process}\label{sec:markov-kernel} Let $p\in (0,1]$. The \textit{BGRW} process with parameter $p$ is a Markov chain with transition kernel $K_p$. To define this kernel, given $(T,x)$, we sample~$(T',x')\sim K_p((T,x),\cdot)$ as follows. 
\begin{enumerate}
\item Let $m=\min\N\backslash V(T)$. With probability $p$, set $V(T') = V(T)\cup\{m\}$ and $E(T')=E(T)\cup\{\{m,x\}\}$. Otherwise, set $T'=T$.
\item Conditionally on the above, let $x'$ be a uniformly chosen neighbor of $x$ in $T'$.
\end{enumerate}
It is easy to see that this does define a valid Markov transition kernel on $\Omega$. We use $(T_t,X_t)_{t\geq 0}$ to denote a trajectory of $K_p$ and $\Pd_{\mu,p}$ and $\Ed_{\mu,p}$ to denote probabilities and expectations when~$(T_0,X_0)\sim \mu$. If $\mu$ is a point mass on $(S_0,x_0)$, we replace $\mu$ by $S_0,x_0$ in the subscripts. 

\begin{remark}Most of the time we will ignore this formal definition of the process and stick with the informal version presented in the introduction. We will later need to define this process on the set of rooted trees up to isomorphisms. See Section \ref{sec:finiteinfinite} for details.\end{remark}

\section{Nontrivial distance from the root}\label{sec:powersoflog}

In this section we show the following property of our process:  for any $y\in T_0$, and any constant $M>0$, the random walker will most likely reach distance $\geq \log^Mn$ from $y$ in at most $n$ steps. 
 
	
\begin{lemma}\label{lem:stretchdistance}For any $M>0$ and  $0<p\leq 1$ there exist $n_0=n_0(p,M)\in\N$, depending only on $p$ and $M$ such that, for all $n\ge n_0$, all finite trees $T_0$ and all $x_0,y\in T_0$,
	\[
	\Pd_{T_0,x_0,p} \left( \exists m\le n, \dist{m}{X_m}{y} \ge \log^Mn\right) \ge 1-e^{-n^{1/4}}.
	\]\end{lemma}

For the proof of this Lemma, we make a definition. 


\begin{definition}\label{def:admissible}Say that $M>0$ is admissible if the conclusion of the Lemma holds for this specific value of $M$. That is, $M$ is admissible if, for any $0<p\leq 1$, there exist $n_0=n_0(p,M)\in\N$,  depending only on $p$ and $M$ such that, for all $n\ge n_0$, all finite trees $T_0$ and all $x_0,y\in T_0$,
	\[
	\Pd_{T_0,x_0,p} \left( \exists m\le n, \dist{m}{X_m}{y} \ge \log^Mn\right) \ge 1-e^{-n^{1/4}}.
\]\end{definition}

Clearly, if $M$ is admissible, so are all $0<M'<M$. Our Lemma follows from the next two Propositions. 


\begin{proposition}[``Easy"; proof in subsection \ref{sec:proof.easypart}]\label{prop:easypart} $M=\frac{1}{2}$ is admissible.\end{proposition}
\begin{proposition}[``Harder"; proof in subsection \ref{sec:proof.harderpart}]\label{prop:harderpart} If $M\geq \frac{1}{2}$ is admissible, so is $M+1/2$.\end{proposition}

In between these two statements, we will also need to prove an intermediate result. It basically says that, when $M$ is admissible and $x_0$ is ``far"~from $y$, then it is likely that the distance from the walker to $y$ will increase by at least one unit by time $n$. This probability is large enough that we are likely to see many such increases in a small time window. 


\begin{claim}[Small growth in distance; proof in subsection \ref{sec:proof.growby1}]\label{claim:growby1} Assume $M\geq \frac{1}{2}$ is admissible (cf. Definition \ref{def:admissible}). Then there exist $n_1(p,M)\in\N$ such that, for $n\geq n_1(p,M)$, the following property holds. Take a finite tree $T$ and $x,y\in T$ with $\dist{}{x}{y}\geq \log^M n$. Then 
\[\Pd_{T,x,p}\left(\dist{t}{X_t}{y} = \dist{}{x}{y}+1\mbox{ for some }t\leq n\right)\geq 1  - 2 \frac{(\log\log n)^2}{\log n}.\]
\end{claim}

As we will see, the ``harder"~Proposition \ref{prop:harderpart} follows from this claim and the assumption that $M$ is admissible applied to time $\approx \sqrt{n}$ instead of $n$. 

Throughout this section we will use the following simple and standard Lemma, which we prove in the Appendix for completeness.  

\begin{lemma}[Proof in Appendix \ref{sec:appendixstochasticdomination}]\label{lem:stochdom}Suppose $(I_j)_{j\in\N\backslash\{0\}}$ are indicator random variables. Assume $\mu$ is such that $\Pd\left(I_1=1\right)\geq \mu$ and
\[\forall j>1\,:\,\Pd\left({I_j=1\mid I_1,\dots,I_{j-1}}\right)\geq \mu.\]
Then for any $k,m\in\N\backslash\{0\}$ 
\[\Pd\left(\mbox{at least $k$ consecutive $1$'s in the sequence $(I_j)_{j=1}^m$}\right)\geq 1 - (1-\mu^k)^{\lfloor m/k\rfloor}.\]\end{lemma}

\subsection{The ``easy" proposition}\label{sec:proof.easypart}

\begin{proof}[Proof of Proposition \ref{prop:easypart}] It suffices to show that there exist $n_0=n_0(p)\in \N$  such that, for all $n\geq n_0$, for  any initial tree $T_0$ and initial vertex $x_0 \in T_0$,
\[\Pd_{T_0,x_0,p}\left(\exists t\leq n\,:\, \dist{t}{X_t}{y}\geq \log^{\frac{1}{2}} n\right)\geq 1 - e^{-n^{1/4}}.\]
To do this, we define, for each time $1\leq t\leq n$, 
\[I_t:=\Ind{}\{\dist{t}{X_t}{y} = \dist{t-1}{X_{t-1}}{y}+1\}.\]
Notice that we have the following inclusion of events
\[
\left\lbrace \exists t\leq n\,:\, \dist{t}{X_t}{y}\geq \log^{\frac{1}{2}} n\right\rbrace \supset \left\lbrace\mbox{at least $\lceil \log^{\frac{1}{2}} n\rceil$ consecutive 1's in $(I_t)_{t=1}^n$}\right \rbrace.
\]
So we will apply Lemma \ref{lem:stochdom}. The point is that, by the Markov property:
\begin{equation}\label{eq:easyboundindicator}
\Pd_{T_0,x_0,p}\left(I_t=1\mid I_1,\dots,I_{t-1}\right)\geq \inf_{T,x}\Pd_{T,x,p}\left(I_1=1\right)\geq \frac{p}{2}.\end{equation}
Indeed, the the value $p/2$ is achieved when $x$ is a leaf of $T$. In that case $I_1=1$ only occurs when a new neighbor is created for $x$ (with probability $p$) and then the walker jumps to that neighbor (with probability $1/2$). If $x$ is not a leaf, then the probability is $>1/2$. 

Combining (\ref{eq:easyboundindicator}) with Lemma \ref{lem:stochdom}, we obtain:
\begin{equation}\label{ineq:c1}
\Pd_{T_0,x_0,p}\left(\exists t\leq n\,:\,\dist{t}{X_t}{y}\geq \log^{\frac{1}{2}} n\right)\geq 1 - \left(1 - \left(\frac{p}{2}\right)^{\lceil \log^{\frac{1}{2}} n\rceil}\right)^{\left\lceil \frac{n}{\lceil  \log^{\frac{1}{2}} n\rceil}\right\rceil}.
\end{equation}
Now, observe that
\begin{equation*}
\begin{split}
	\left(1 - \left(\frac{p}{2}\right)^{\lceil \log^{\frac{1}{2}} n\rceil}\right)^{\left\lceil \frac{n}{\lceil  \log^{\frac{1}{2}} n\rceil}\right\rceil} & \le \exp \left \lbrace - \exp \left \lbrace \log n\left(1-\frac{\log\left(\frac{2}{p}\right)\left(\sqrt{\log n} + 1\right)}{\log n} - \frac{\log \log n}{\log n} \right)\right \rbrace \right \rbrace.
\end{split}
\end{equation*}
Choosing $n_0$ sufficiently large such that $ \frac{\log\left(\frac{2}{p}\right)\left(\sqrt{\log n} + 1\right)}{\log n} + \frac{\log \log n}{\log n} < \frac{3}{4}$,  the lower bound (\ref{ineq:c1}) is $\geq 1-e^{-n^{1/4}}$ for all $n\geq n_0$. For later purposes we point out that, 
$n_0(p)$ is not-increasing as a function of $p$.  
\end{proof}

\subsection{A key claim}\label{sec:proof.growby1}

We now come to the proof of Claim \ref{claim:growby1}, which connects the ``easy"~and ``hard"~proposition in the proof. 

\begin{proof}[Proof of Claim \ref{claim:growby1}]Consider the vertex $y_*$ on the unique path from $x$ to $y$ with~$\dist{}{x}{y_*}=\lceil \log^M n\rceil -1$. We recall that the admissibility of $M$ implies, for all $0<p\leq 1$, the existence of $n_0$ depending only on $M$ and $p$ such that:
\begin{equation}\label{eq:admissibleystar}\forall n\geq n_0\,:\,\Pd_{T,x,p}(\exists t\leq n\,:\, \dist{t}{X_t}{y_*} \geq \log^M n) \geq 1 - e^{-n^{1/4}}.
\end{equation}
Now let $F$ (for failure) denote the event that $\dist{t}{X_t}{y_*} \leq \dist{}{x}{y}$ for all $t\leq n$. Let~$\tau_{y_*}$ the hitting time of $y_*$
\[\tau_{y_*}:=\inf\{t\in\N\,:\, X_t=y_*\}\in \N\cup\{+\infty\}.\]
Define $\tau^{+k}_x$ the $k$-th return time to $x$. That is, we set $\tau^{+0}_x:=0$, and then (recursively for $k\in\N\backslash\{0\}$):
\[\tau^{+k}_x:=\left\{\begin{array}{ll}+\infty, & \mbox{ if }\tau_x^{+(k-1)}=+\infty; \\ \inf\{t> \tau^{+(k-1)}_x\,:\,X_t=x\}\in \N\cup\{+\infty\}, & \mbox{otherwise}.\end{array}\right.\]
Note that, for any $k$, we may upper bound the probability of $F$ by three terms, which we will bound separately.
\begin{equation}\label{eq:failureclaim3}
\Pd_{T,x,p}(F)\leq\Pd_{T,x,p}(F\cap \{\tau_{y_*}>n\}) + \Pd_{T,x,p}(F\cap \{\tau^{+k}_x<\tau_{y_*}\leq n\}) + \Pd_{T,x,p}(\tau^{+k}_x>\tau_{y_*}).
\end{equation}

We start with the first term in the RHS. Note that if $\tau_{y_*}>n$, then the unique path from $X_t$ to $y$ passes through $y_*$ at all times $t\leq n$. Therefore, 
\[
\forall t\leq n\,:\,\dist{t}{X_t}{y} - \dist{}{x}{y} = \dist{t}{X_t}{y_*} - \dist{}{x}{y_*}.
\]
In particular, when $F\cap \{\tau_{y_*}>n\}$ holds, we must have
\[\forall t\leq n\,:\, \dist{t}{X_t}{y_*} \leq  \dist{}{x}{y_*} <\log^M n \mbox{, i.e. the event in (\ref{eq:admissibleystar}) does not hold}.\]
We conclude:
\begin{equation}\label{eq:ystarnotcrossed}
\Pd_{T,x,p}(F\cap \{\tau_{y_*}>n\})\leq e^{-n^{1/4}}.
\end{equation}
We now consider the second term in the RHS of (\ref{eq:failureclaim3}). In order for $F\cap \{\tau^{+k}_x<\tau_{y_*}\leq n\}$ to take place, it must be that $X_t$ returns at least $k$ times to $x$ before visiting $y_*$ {\em but} never gets to jump to a neighbor $x'$ of $x$ with ${\rm dist}(x',y)={\rm dist}(x,y)+1$. Now, at each return to $x$, the probability that $X_t$ jumps to such a neighbor {\em conditionally on the process up to that point} is at least $p/2$: the probability of creating a leaf and then not jumping in the direction of $y$\footnote{This is similar to what we did when we proved equation (\ref{eq:easyboundindicator}) in the proof of Proposition \ref{prop:easypart}.}. We deduce:
\begin{equation}\label{eq:Fwithmanyreturns}
\Pd_{T,x,p}(F\cap \{\tau^{+k}_x<\tau_{y_*}\leq n\})\leq \left(1- \frac{p}{2}\right)^k.\end{equation}

Finally, we come to the third term in the RHS of (\ref{eq:failureclaim3}). Consider the walk $X_t$ at the sequence of time steps that it spends on  the path from $x$ to $y$, up to the time $\tau_{y_*}$. The resulting process is a simple random walk on the path with potential delays and reflecting barriers: indeed, since our graph is a tree at all times, whenever the random walk leaves the path, it must return to it (if it does return) at the same point that was last visited. Now, in order that~$\tau_x^{+k}>\tau_{y_*}$, it must be that this simple random walk on the path hits $y_*$ before returning $k$ times to $x$. Since the path has length $\dist{}{x}{y}=\lceil \log^Mn\rceil-1$, the probability of this happening is $\leq k/ (\lceil \log^Mn\rceil-1)$, and we obtain:
\begin{equation}\label{eq:Fwithfewreturns}\Pd_{T,x,p}(\tau^{+k}_x>\tau_{y_*})\leq \frac{k}{\lceil \log^Mn\rceil-1}.\end{equation}

Combining the terms in (\ref{eq:ystarnotcrossed}), (\ref{eq:Fwithmanyreturns}) and (\ref{eq:Fwithfewreturns}), we obtain a bound in (\ref{eq:failureclaim3}):
\begin{equation}\label{ineq:tomjobim}
\begin{split}
\Pd_{T,x,p}(F) & \leq e^{-n^{1/4}} + \left(1- \frac{p}{2}\right)^k + \frac{k}{\lceil \log^Mn\rceil-1} \\
& \leq \frac{\lceil(\log \log n)^2\rceil}{\log^Mn} \left[ e^{-n^{1/4}}\frac{\log^Mn}{\lceil(\log \log n)^2\rceil} + e^{-\frac{p}{2}(\log \log n)^2}\frac{\log^Mn}{\lceil(\log \log n)^2\rceil} + \frac{\log^Mn}{\log^Mn -1}\right],\\
\end{split}
\end{equation}
where the last inequality follows by choosing $k:=\lceil (\log\log n)^2\rceil$. Choosing $n_1 = n_1(p,M)$ sufficiently large in (\ref{ineq:tomjobim}) such that 
$ e^{-n^{1/4}}\frac{\log^Mn}{\lceil(\log \log n)^2\rceil} + e^{-\frac{p}{2}(\log \log n)^2}\frac{\log^Mn}{\lceil(\log \log n)^2\rceil} + \frac{\log^Mn}{\log^Mn -1}<2$ we obtain that, for all $n\geq n_1$,
\begin{equation}
\begin{split}
\Pd_{T,x,p}(F) & \leq  2 \frac{(\log \log n)^2}{\log^Mn}.
\end{split}
\end{equation}
 This proves the claim. Note that $n_1(p,M)$ is not-increasing in $p$.
\end{proof}


\subsection{The ``harder" proposition}\label{sec:proof.harderpart}

We now prove Proposition \ref{prop:harderpart}, thus finishing the proof of Lemma \ref{lem:stretchdistance}.

\begin{proof}[Proof of Proposition \ref{prop:harderpart}] Fix an admissible $M\geq \frac{1}{2}$. Define (with hindsight) two time lengths:
\begin{equation}\label{eq:timelengths2}t_*:=\lceil \log^{M+1/2}\,n\rceil,\, \ell_*:=\left\lfloor\sqrt{n}\right\rfloor.\end{equation}
Both of these time lengths depend on $n$, but we omit this dependency from the notation to avoid clutter. The definitions of these quantities will probably seem misterious, but we comment on them in due time; see Remarks \ref{rem:lstarlarge}, \ref{rem:lstarsmall} and \ref{rem:tstarnottoolarge} below. 

The fact that $M$ is admissible and $\ell_*\to +\infty$ when $n\to +\infty$ implies that, for any $0<p\leq 1$, there exists $n_2(p,M)$ such that if $n\geq n_2(p,M)$, we have the following properties for any finite tree $T_0$ and any $x_0\in T_0$. 
\begin{enumerate}
\item By the definition of admissibility (Definition \ref{def:admissible}) applied with $\ell_*$ replacing $n$: 
\begin{equation}\label{eq:Madmissiblegives}
	\Pd_{T_0,x_0,p} \left( \exists m \le \ell_*, \dist{m}{X_m}{y} \ge 2^{-M}\log^{M}n \right) \ge 1 - e^{-\ell_*^{1/4}}\geq 1-e^{-\frac{1}{2}n^{1/8}}
\end{equation}
\item By Claim \ref{claim:growby1} applied with $\ell_*$ replacing $n$, if $T$ is finite and $x,y\in T$~with $\dist{}{x}{y}\geq 2^{-M}\log^M n$, then: 
\begin{equation}\label{eq:Claimgives}\Pd_{T,x,p}\left(\dist{t}{X_t}{y} = \dist{}{x}{y}+1\mbox{ for some }t\leq \ell_*\right)\geq 1  - 2\frac{(\log\log n)^2}{2^{-M}\log^M n}.\end{equation}
Again we used that $\log \ell_*\approx \log n$.
\end{enumerate}

\begin{remark}[$\ell_*$ is large enough]\label{rem:lstarlarge}In the above we used the fact that $\log\ell_*\approx \log n$ to guarantee that the two events under consideration have high probability. We will later need that $n/\ell_*t_*$ is large; see Remark \ref{rem:lstarsmall} below.\end{remark}

We now define a sequence of stopping times $\sigma_j$ and indicators $I_j$. Intuitively, we will want that $\dist{\sigma_j}{X_{\sigma_j}}{y}>\dist{\sigma_{j-1}}{X_{\sigma_{j-1}}}{y}$, and we will signal such a success by setting $I_j=1$. More formally, for $j=0$, $\sigma_j=0$ and $I_j=0$. We define $\sigma_j$ and $I_j$ for $j>0$, with the following two choices.
\begin{enumerate} 
\item[{\bf (a)}] {\em When $X_{\sigma_{j-1}}$ is too close to $y$:} if $\dist{\sigma_{j-1}}{X_{\sigma_{j-1}}}{y}<2^{-M}\log^M n$, we let:
\begin{eqnarray}\label{eq:defsigmajclose}\sigma_j&:=&\inf\{t\geq \sigma_{j-1}\,:\,\dist{t}{X_{t}}{y} \geq 2^{-M} \log^{M}n\mbox{ or }t-\sigma_{j-1} = \ell_*\}\;\\  \label{eq:defIjclose}
I_j &=& \Ind{}\{\dist{\sigma_j}{X_{\sigma_j}}{y} \geq 2^{-M}\log^{M}n\}
\end{eqnarray}
\item[{\bf (b)}] {\em When $X_{\sigma_{j-1}}$ is not too close to $y$:} if $\dist{\sigma_{j-1}}{X_{\sigma_{j-1}}}{y}\geq 2^{-M} \log^M n$, we set
\begin{eqnarray}\label{eq:defsigmajfar}\sigma_j&:=&\inf\{t\geq \sigma_{j-1}\,:\,\dist{t}{X_{t}}{y} =\dist{\sigma_{j-1}}{X_{\sigma_{j-1}}}{y} + 1\mbox{ or }t-\sigma_{j-1} = \ell_*\}\;\\  \label{eq:defIjfar}
I_j &=& \Ind{}\{\dist{t}{X_{t}}{y} =\dist{\sigma_{j-1}}{X_{\sigma_{j-1}}}{y} + 1 \}\end{eqnarray}\end{enumerate}

We note some basic properties of these new random variables. First, we have the deterministic bound $\sigma_j-\sigma_{j-1}\leq \ell_*$. In particular, $\sigma_j\leq \ell_*\,j$ for all $j$.
Second, we observe that, if we have $m+1$ consecutive successes, i.e. if $I_{j_0} = I_{j_0+1} = \dots = I_{j_0 + m} = 1$ for some $j_0,m\in\N\backslash\{0\}$, then:
\[
\dist{\sigma_{j_0}}{X_{\sigma_{j_0}}}{y}\geq 2^{-M}\log^{M}n
\]
and
\[
\dist{\sigma_{j_0+m}}{X_{\sigma_{j_0+m}}}{y}= \dist{\sigma_{j_0+m-1}}{X_{\sigma_{j_0+m-1}}}{y}+1=\dots = \dist{\sigma_{j_0}}{X_{\sigma_{j_0}}}{y} + m.
\]
We thus arrive at the following crucial observation (recall the time lengths $t_*,\ell_*$ in (\ref{eq:timelengths2})).

\begin{obs} \label{obs:consecutive} Assume that there are $t_*+1$ consecutive ones in the sequence $(I_j)_{j=1}^{\lfloor n /\ell_*\rfloor}$. Then, there exists $t\leq n$ such that
\[ \dist{t}{X_t}{y}\geq t_*\geq \log^{M+1/2}n.\]\end{obs}

Indeed, if we have $t_*+1$ consecutive ones in this sequence there exists $j_0$ with $(j_0+t_*)\ell_*\leq n$ and $I_{j_0} = I_{j_0+1} = \dots = I_{j_0 + t_*}=1$. What we are after is to show that the probability of $t_*+1$ consecutive ones is very likely.  The upshot is that we may apply Lemma \ref{lem:stochdom} above to control the probability that the RW reaches distance $\log^{M+1/2}n$ from $y$ in $n$ time steps, i.e. that $M+1/2$ is admissible. 
\begin{remark}\label{rem:lstarsmall}Note that for the Lemma to be effective we need 
$n/\ell_*\gg 1$ so that there are ``many indicators"~to consider, and also that $n/\ell_*t_*\gg 1$. It will later become clear that (for $t_*$ polylogarithmic) it suffices that $n/\ell_*\geq n^{1/4+c+o(1)}$ for some $c>0$.\end{remark}

To apply Lemma \ref{lem:stochdom}, we must lower bound $\Pd_{T_0,x_0,p}(I_1=1)$ and $\Pd_{T_0,x_0,p}\left(I_j=1\mid I_1,\dots,I_{j-1}\right)$. By the strong Markov property, these quantities are lower bounded by
\[\inf_{T,x}\Pd_{T,x,p}\left(I_1=1\right)\]
Crucially, conditions {\bf (a)} and {\bf (b)} in the definition of $\sigma_1$ and $I_1$ correspond precisely to the situations in the two bounds (\ref{eq:Madmissiblegives}) and (\ref{eq:Claimgives}) (respectively). It follows that, for any $0<p\leq 1$,  there exists a $n_3=n_3(p,M)$ such that,  for all $n\geq n_3$, 
\[\inf_{T,x}\Pd_{T,x,p}\left(I_1=1\right)\geq \mu:= 1 - 2\,\frac{(\log\log n)^2}{2^{-M}\log^Mn}.\]

Lemma \ref{lem:stochdom} implies that:
\[\Pd_{T_0,x_0,p}\left(\mbox{at least $t_*+1$ consecutive ones in $(I_j)_{j=1}^{\lfloor n /\ell_*\rfloor}$}\right)\geq 1  - (1-\mu^{t_*+1})^{\left\lfloor \frac{\lfloor n/\ell_*\rfloor}{t_*+1}\right\rfloor}.\]

What is left to show is that, for any $0<p\leq 1$, there exists $n_4=n_4(p,M)$ such that 
\begin{equation}\label{eq:jorge-bem}\forall n\geq n_4\,:\,
(1-\mu^{t_*+1})^{\left\lfloor \frac{\lfloor n/\ell_*\rfloor}{t_*+1}\right\rfloor}\leq  \exp(-n^{1/4}) ,
\end{equation}
which would imply 
\[\forall n\geq n_4\,:\,\Pd_{T_0,x_0,p}\left(\mbox{at least $t_*$ consecutive ones in $(I_j)_{j=1}^{\lfloor n /\ell_*\rfloor}$}\right)\geq 1- \exp(-n^{1/4}).\]
The latter bound, which is uniform in $T_0,x_0$, together with
Observation \ref{obs:consecutive},  will  prove that $M+1/2$ is admissible.   
\begin{remark}\label{rem:tstarnottoolarge}The important thing here is that, because $t_*=\log^{M+1/2}n $, the probability of $t_*+1$ consecutive ones, i.e. $\mu^{t_*+1}$, goes to $0$ slowly. Remark \ref{rem:lstarsmall} shows we are considering ``polynomially many indicators", so such not-too-small probabilities make a long run of ones very likely. 
\end{remark}

Let us show that for our choices of $t_*$ and $\ell_*$ we can find an $n_4$ sufficiently large such that Equation~\eqref{eq:jorge-bem}  holds. We point out that we will implicitly assume $n_4> n_3$ so to guarantee that
$\inf_{T,x}\Pd_{T,x,p}\left(I_1=1\right)\geq \mu$.
It is for this reason that $n_4$ will depend on $p$ (in particular, $n_4$ will be non-increasing in $p$).
\[
(1-\mu^{t_*+1})^{\left\lfloor \frac{\lfloor n/\ell_*\rfloor}{t_*+1}\right\rfloor}
\leq \exp\left(- \mu^{t_*+1} \frac{n^{1/2}}{2(\log^{M+1/2} n + 1)}\right) 
\]

Since $t_*=\lceil \log^{M+1/2} n\rceil$, for sufficiently large $n$ we have
\[ \mu^{t_*+1} \geq \exp\left( -\left( 2(\log\log n)^2 2^{M+1} + o((\log\log n)^2)\right) \sqrt{\log n}\right)\mu  
=n^{-o(1)}\mu\,
\]
where we are using that $\left(1-\frac{b_n}{a_n}\right)^{a_n}\approx e^{-b_n -o(b_n)}$ for sufficiently large $n$ whenever $a_n, b_n \to + \infty$ and $b_n=o(a_n)$. Thus,
\[
(1-\mu^{t_*+1})^{\left\lfloor \frac{\lfloor n/\ell_*\rfloor}{t_*+1}\right\rfloor}\leq \exp\left(-n^{1/4}\, \frac{n^{1/4-o(1)}\mu}{2(\log^{M+1/2} n+1)}\right) 
\leq \exp\left(-n^{1/4}\right)
\]
for large enough $n$.


\end{proof}

\section{The loop process, or why it's hard to go back}\label{sec:loopprocess}

Now we stop the discussion about the \textit{BGRW} process, to introduce a simpler process which will help us to understand how long the walker $X$ stays on specific subgraphs of the random trees $\{T_n\}_{n \in \N}$. Roughly speaking, a loop process on an initial graph $G$ is a random walker such that at each step adds a loop to its position according to a coin and then chooses uniformly one edge of its current position to walk on. In other words, the process is quite similar to the \textit{BGRW} but here the walker adds loops instead of leaves, which makes it possible for it to stand still.

\begin{figure}[h]
	\centering
	\includegraphics[width=0.95\linewidth]{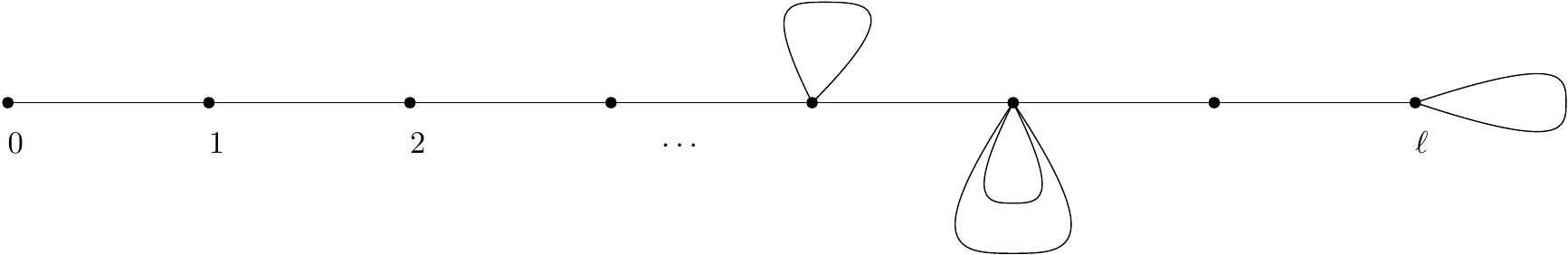}
	\caption[A backbone]{A \textit{backbone} of length $\ell$}
	\label{fig:backbone}
\end{figure}

We will be particularly interested in the loop process over specific graphs which we define before the definition of the process itself. We call a finite graph $\mathcal{B}$ a \textit{backbone} of length $\ell$ if~$\mathcal{B}$ is a path of length $\ell$ having a loop attached to its $\ell + 1$-th vertex and possibly to its other vertices, see Figure~\ref{fig:backbone}.

In this section, we will abuse the graph terminology saying \textit{degree} of a vertex even though we do not count loops twice. We then reserve the special notation $\text{deg}_t(i)$ to denote the number of edges attached to vertex $i$ at time $t$. 

The model has one parameter $p \in [0,1]$ which is the parameter of a sequence $\{Z_t\}_{t\in \N}$ of \textit{i.i.d} random variables with law $\text{Ber}(p)$. The loop process on a backbone of length~$\ell$ is a Markov chain $\{( \mathcal{B}_t ,X^{\text{loop}}_t)\}_{t \in \mathbb{N}}$ where $\mathcal{B}_t$ denotes the resultant backbone at time $t$ and~$X^{\text{loop}}_t$ is one of its~$\ell+1$ vertices.  The loop process is also defined inductively according to the update rules below
\begin{enumerate}
	\item Obtain $\mathcal{B}_{t+1}$ from $\mathcal{B}_t$ by adding a new loop to $X^{\text{loop}}_t$ whether $Z_{t+1} = 1$;
	\item Choose uniformly one edge attached to $X^{\text{loop}}_t$ in $\mathcal{B}_{t+1}$. Whether the chosen edge is a loop set $X^{\text{loop}}_{t+1}$ as $X^{\text{loop}}_t$. Otherwise, $X^{\text{loop}}_{t+1}$ becomes the $X^{\text{loop}}_t$ neighbor.
\end{enumerate}
We stress out to the index $t+1$ of $\mathcal{B}$ on~(2). It means that we may add a new loop at (1) and then choose it at (2).

For a fixed backbone $\mathcal{B}$ of length $\ell$ and $i \in \{0,1,\cdots, \ell\}$ we let $\Pd_i(\cdot)$ denote the law of~$\{( \mathcal{B}_t ,X^{\text{loop}}_t)\}_{t \in \mathbb{N}}$ when~$( \mathcal{B}_0 ,X^{\text{loop}}_0) \equiv (\mathcal{B},i)$.

Once we have defined the loop process, we are interested on the time it takes to go from one end of the backbone to another. More precisely, we would like to obtain bounds on the stopping time bellow
\begin{equation}\label{def:etaloop0}
\eta^{\text{loop}}_0 := \inf\left \lbrace t \ge 1 \middle |X^{\text{loop}}_t = 0 \right \rbrace
\end{equation}
when the process starts from $\ell$. The next Lemma gives us some estimates.
\begin{lemma}\label{lemma:xlooptime} There exist positive constants $c_1$ and $c_2$ depending on $p$ only, such that, for all integer $K$
	\[
	\Pd_{\ell} \left( \eta^{\text{loop}}_0 \le e^{K}\right) \le 1 - \left(1-\frac{1}{2\ell}\right)^{c_2K} + e^{-c_1K}
	\]
\end{lemma}
\begin{proof} We need some notation and definitions. Consider the following sequence of stopping times
	\begin{equation}
	\begin{split}
	& \tau_0 \equiv 0, \\
	& \tau_k := \inf\left \lbrace t > \tau_{k-1} | X^{\text{loop}}_t \neq X^{\text{loop}}_{\tau_{k-1}}\right \rbrace.
	\end{split}
	\end{equation}
	Observe that the probability of $X^{\text{loop}}_t$ leaving its current position is at least $1/t$ since the degree of a vertex at time $t$ is at most $t$. This implies that $\tau_k$ is finite a.s. for all $k$. This allow us to define the process
	\begin{equation}
	Y_k := X^{\text{loop}}_{\tau_k}.
	\end{equation}
	Note that by strong Markov Property, $\{Y_k\}_k$ is a simple \textit{RW} on $\{0,1,\cdots,\ell\}$ with reflecting barriers. 
	Regarding the process $\{Y_k\}_k$ we let $\s$ be the following stopping time
	\begin{equation}
	\s :=  \inf\left \lbrace k > 0 | Y_k = 0\right \rbrace.
	\end{equation}
	Observe that $\eta^{\text{loop}}_0 = \tau_{\s}$.
	
	We prove the Lemma by showing that $X^{\text{loop}}$ spends at least $\exp\{K\}$ steps on the vertex~$\ell$. More precisely, we prove that the degree of $\ell$ at time $\tau_{\s}$ is at least $\exp\{K\}$, \textit{w.h.p.} To do this, first observe that the degree of a vertex may be written in terms of $Y_k$ and $\D \tau_k$ as follows
	\begin{equation}
	\mathrm{deg}_{\tau_n}(\ell) = 2+ \sum_{k=0}^{n} \mathbb{1}\{Y_k = \ell\}\left(\sum_{m=\tau_k}^{\tau_{k+1}-1}Z_m\right) = 2+ \sum_{k=0}^{n} \mathbb{1}\{Y_k = \ell\}\mathrm{Bin}\left(\D \tau_k, p \right).
	\end{equation}
	Also notice that if $Y_k = \ell$, then the number of steps $X^{\text{loop}}$ spends on $\ell$ is exactly $\D \tau_k$, which satisfies
	\[
	\D \tau_k \ge \mathrm{Geo}\left(1/\mathrm{deg}_{\tau_k}(\ell)\right).
	\]
	To see the bound above, consider the random variable which counts the number of steps $X^{\text{loop}}$ spends on $\ell$ by choosing only the loops which were already attached on $\ell$ when $X^{\text{loop}}$ arrived at $\ell$. This random variable is clearly smaller than $\D \tau_k$ and follows a geometric distribution of parameter $1/\mathrm{deg}_{\tau_k}(\ell)$. Thus, we have
	\begin{equation}
	\mathrm{deg}_{\tau_n}(\ell) \ge \sum_{k=0}^{n} \mathbb{1}\{Y_k = \ell\}\mathrm{Bin}\left(\mathrm{Geo}\left(1/\mathrm{deg}_{\tau_k}(\ell)\right), p \right).
	\end{equation}
	Regarding the random variables $\mathrm{Bin}\left(\mathrm{Geo}\left(1/\mathrm{deg}_{\tau_k}(\ell)\right), p \right)$, we claim that
	\begin{claim} For all $\e \in (0,1)$, there exists a positive constant $q=q(\e,p)$ such that
		\begin{equation}
		\Pd_{\ell}\left( Y_k = \ell, \mathrm{Bin}\left(\mathrm{Geo}\left(1/\mathrm{deg}_{\tau_k}(\ell)\right), p \right) \ge \e p\mathrm{deg}_{\tau_k}(\ell) \middle | \mathcal{F}_{\tau_k}\right) \ge q\mathbb{1}\{Y_k = \ell\}.
		\end{equation}
	\end{claim}
	\begin{claimproof} To simplify our writing, write
		\[
		G_k := \mathrm{Geo}\left(1/\mathrm{deg}_{\tau_k}(\ell)\right); \; d_k := \mathrm{deg}_{\tau_k}(\ell)
		\]
		and let $\tilde{\mathcal{F}}_{\tau_k}$ be the $\s$-algebra generated by $G_k$ and $\mathcal{F}_{\tau_k}$. Now, by Chernoff bounds, we have that
		\begin{equation}
		\begin{split}
		\Pd_{\ell}\left( Y_k = \ell, \mathrm{Bin}\left(G_k, p \right) \ge \e pd_k, G_k \ge  d_k \middle | \tilde{\mathcal{F}}_{\tau_k}\right) &\ge \left(1- e^{-(1-\e)^2pG_k}\right)\mathbb{1}\{Y_k = \ell, G_k \ge d_k \} \\
		& \ge \left(1- e^{-(1-\e)^2pd_k}\right)\mathbb{1}\{Y_k = \ell, G_k \ge d_k \}.
		\end{split}
		\end{equation}
		Recall that $d_k$ is greater than $2$ for all $k$. So, taking the conditional expectation \textit{wrt} $\mathcal{F}_{\tau_k}$ on the above inequality yields
		\begin{equation}
		\begin{split}
		\Pd_{\ell}\left( Y_k = \ell, \mathrm{Bin}\left(G_k, p \right) \ge \e pd_k, G_k \ge  d_k \middle | \mathcal{F}_{\tau_k}\right)
		& \ge \left(1- e^{-2(1-\e)^2p}\right)\Pd_{\ell} \left( Y_k = \ell, G_k \ge d_{\tau_k} \middle | \mathcal{F}_{\tau_k} \right) \\
		& \ge \left(1- e^{-2(1-\e)^2p}\right)\left(1- e^{-1}\right)\mathbb{1}\{Y_k = \ell\},
		\end{split}
		\end{equation}
		which proves the claim.
	\end{claimproof}
	
	The above claim tells us that conditionally on the past, every time $\{Y_k\}_k$ visits $\ell$ it has probability at least $q$ of increasing the degree of $\ell$ by a factor of at least $1+\e p$. This points out that the degree of $\ell$ must be at least exponential of the number of visits $\ell$ receives from $\{Y_k\}_k$. So, let $N_{\s}(\ell)$ the number of visits made by $Y$ to $\ell$ before it reaches vertex~$0$. Since $Y$ is a simple random walk on $\{0,1,\cdots,\ell\}$, $N_{\s}(\ell)$ follows a geometric distribution of parameter~$1/2{\ell}$. Moreover, the random variable $W$ that counts how many times we have successfully multiplied the degree of $\ell$ by $1+\e p$ may be written as follows
	\begin{equation}
	W := \sum_{k=0}^{\s}\mathbb{1}\{Y_k = \ell\}\mathbb{1}\{\mathrm{Bin}\left(\mathrm{Geo}\left(1/\mathrm{deg}_{\tau_k}(\ell)\right), p \right) \ge \e p\mathrm{deg}_{\tau_k}(\ell)\}
	\end{equation}
	and dominates a random variable distributed as $\text{Bin}(N_{\s}(l),q)$. Consequently, by Chernoff bounds
	\begin{equation}\label{ineq:w}
	\begin{split}
	\Pd_{\ell}\left(W \le \frac{K}{\log(1+\e p)}\right) & \le \Pd_{\ell}\left(\mathrm{Bin}(N_{\s}(\ell),q) \le \frac{K}{\log(1+\e p)}\right) \\
	& \le \Pd_{\ell}\left(\mathrm{Bin}(N_{\s}(\ell),q) \le \frac{K}{\log(1+\e p)}\middle |N_{\s}(\ell) \ge \frac{2K}{q\log(1+\e p)} \right) \\
	& + \Pd_{\ell}\left(N_{\s}(\ell) < 2q^{-1}K/\log(1+\e p)\right) \\
	& \le \exp\{-c_1K\} + 1-\left(1-\frac{1}{2\ell}\right)^{c_2K}.
	\end{split}
	\end{equation}
	Finally, observe that if $W \ge K/\log(1+\e p)$, then $\mathrm{deg}_{\tau_{\s}}(\ell) \ge 2e^{K}$ which implies that $\tau_{\s}$ is at least this amount, finishing the proof.
\end{proof}
The following special case of the above Lemma will be particularly useful to our proposes.
\begin{corollary}\label{cor:xlooptime} There exists a positive constant $C$ depending on $p$ only, such that
	\[
	\Pd_{\ell} \left( \eta^{\text{loop}}_0 \le e^{\sqrt{\ell}}\right) \le \frac{C}{\sqrt{\ell}}.
	\]
\end{corollary}
\begin{proof}Letting $K = \sqrt{\ell}$ on (\ref{ineq:w}) of Lemma~\ref{lemma:xlooptime} we obtain
	\begin{equation}
	\begin{split}
	\Pd_{\ell}\left(W \le \frac{\sqrt{\ell}}{\log(1+\e p)}\right) & \le \Pd_{\ell}\left(\text{ Bin}(N_{\s}(\ell),q) \le \frac{\sqrt{\ell}}{\log(1+\e p)}\right) \\
	& \le \Pd_{\ell}\left(\mathrm{Bin}(N_{\s}(\ell),q) \le \frac{\sqrt{\ell}}{\log(1+\e p)}\middle |N_{\s}(\ell) \ge \frac{2\sqrt{\ell}}{q\log(1+\e p)} \right) \\
	& + \Pd_{\ell}\left(N_{\s}(\ell) \le 2q^{-1}\sqrt{\ell}/\log(1+\e p)\right) \\
	& \le \exp\{-c_1\sqrt{\ell}\} + 1-\left(1-\frac{1}{2\ell}\right)^{c_2\sqrt{\ell}}.
	\end{split}
	\end{equation}
	Since $\left(1-\frac{1}{2\ell}\right)^{c_2\sqrt{\ell}} \approx e^{-c_3/\sqrt{\ell}}$ and $1-e^{-x} \le x$, choosing properly the constants we obtain the desired result.
\end{proof}
\subsection{Coupling the BGRW and the loop process}\label{sec:couplingbgrwloop}
In this subsection we construct a coupling of the \textit{BGRW} and loop processes in such way that the loop process is always closer to the root than the walker $X$. For this, let $T$ be a rooted locally finite tree of height at least~$2(\ell+1)$,~$x$ a vertex such that~$\dist{}{x}{\text{root}} \ge \ell$ and $\degr{}{x} \ge 2$ and $y$ its ancestor at distance $\ell$. Since~$T$ is a tree, there exists only one path $\mathcal{P}$ connecting $x$ to the $y$. With this in mind, we define a graph operation $\mathcal{B}$ which associates to each pair~$(T,x)$ and ancestor $y$ satisfying the aforementioned conditions a backbone $\mathcal{B}(T,y,x)$ of length $\ell$ as follows:
\begin{enumerate}
	\item delete all vertices of $T$ whose distance from $\mathcal{P}$ is at least $2$;
    \item replace each edge $xy\in E(T)$ with $x\in \mathcal{P}$ and $y\not\in\mathcal{P}$ with a loop edge at $x$ (so each edge stemming out of the path becomes a loop).

	\item label the vertices on $\mathcal{P}$ by their distance from $y$ (so $y$ gets label $0$, its neighbor on $\mathcal{P}$ gets label $1$ and so on).
\end{enumerate}
The figure below gives a concrete example of $\mathcal{B}$ in action when $y$ is taken as the root:
\begin{figure}[h]
	\centering
	\includegraphics[height=0.4\textheight]{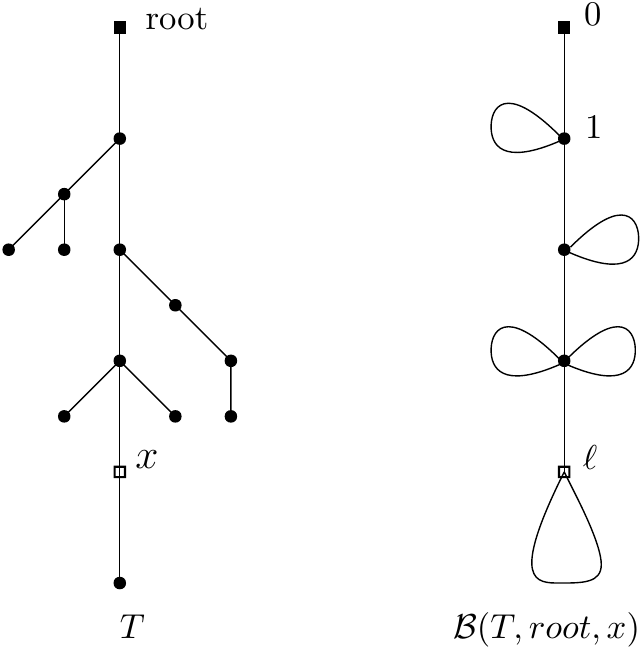}
	\caption{Example of a pair $(T,x)$ transformed into a backbone}
	\label{fig:backbonetransform}
\end{figure}

So far, we have shown that the \textit{BGRW} is capable of reaching long distances -- powers of $\log n$ -- away from the root. Now, we would like to argue that once it has gone so far, it takes too long to return. More generally, if the \textit{BGRW} starts on $T_0,x_0$ and $y$ is one of its ancestors, we would like to obtain lower bounds on the following stopping time
\begin{equation}
\eta_{y} := \inf \left \lbrace n \ge 1 \middle |X_n = y \right \rbrace.
\end{equation}
The way we bound $\eta_{y}$ from below is comparing it with $\eta^{\text{loop}}_0$, which we recall its definition
\[
\eta^{\text{loop}}_0 := \inf  \left \lbrace t \ge 1 \middle |X^{\text{loop}}_t = 0 \right \rbrace.
\]
The next proposition tells us that we may couple the \textit{BGRW} and the loop process in such way that $\eta_{y}$ is greater than $\eta^{\text{loop}}_0$ almost surely.
\begin{proposition}[Coupling]\label{prop:coupling} Let $T_0$ be a rooted locally finite tree, $x_0$ one of its vertices different from the root and $y$ an ancestor of $x_0$. Then, there exists a coupling of $\{(T_n,X_n)\}_{n\in\N}$ starting from~$(T_0,x_0)$ and a \textit{loop} process $\{(\mathcal{B}_n,X^{\text{loop}}_n)\}_{n\in\N}$ starting from~$(\mathcal{B}(T_0,y,x_0),x_0)$ such that
	\[
	\Pd\left(\eta_{y}\ge \eta^{\text{loop}}_0\right) = 1.
	\]
\end{proposition}
\begin{proof}Let $\mathcal{P}$ denote the path connecting $x_0$ to its ancestor $y$ on~$T_0$ and $\ell$ its length. Also consider the following sequence of stopping times,
	\begin{equation}
	\zeta_{0} \equiv 0; \quad \zeta_k := \inf \left \lbrace m > \zeta_{k-1} \middle | X_m \in \mathcal{P}\right \rbrace,
	\end{equation}
	and let $\{W_k\}_{k\in \N}$ be a sequence of \textit{i.i.d.} r.v's independent of the process \textit{BGRW}, such that~$W_k\sim \text{Ber}(p)$. We couple a loop process $\{(\mathcal{B}_k,X^{\text{loop}}_k)\}_{k\in \N}$ to the \textit{BGRW} inductively as follows. Start by defining
	\[
	(\mathcal{B}_0, X^{\text{loop}}_0) := (\mathcal{B}(T_0,y,x_0), x_0)
	\]
	and assume we have defined $\{(\mathcal{B}_i,X^{\text{loop}}_i)\}_{i=0}^{k-1}$ in such way that this random vector is distributed as $k-1$ steps of a loop process starting from $(\mathcal{B}_0, \ell)$. Now, we define $(\mathcal{B}_k,X^{\text{loop}}_k)$ 
	\begin{itemize}
		\item If $\zeta_{k-1} < \infty$, we set $\mathcal{B}_k = \mathcal{B}(T_{\zeta_{k-1}+1},y,x_0)$;
		\item If $\zeta_{k-1} = \infty$, we make use of our independent sequence $\{W_k\}_{k\in \N}$. We modify $\mathcal{B}_{k-1}$ adding a loop to $X^{\text{loop}}_{k-1}$ whether $W_{k} = 1$;
	\end{itemize}
	Observe that if $\zeta_{k-1}$ is finite, then $X_{k-1}^{\text{loop}} = X_{\zeta_{k-1}}$. In fact, when $X_{k-1}^{\text{loop}} = X_{\zeta_{k-2}}$ we know that $X_{\zeta_{k-2}}$ has jumped outside $\mathcal{P}$. Thus, the only way of $X$ coming back to $\mathcal{P}$ is through $X_{\zeta_{k-2}}$, which implies $X_{\zeta_{k-1}} = X_{\zeta_{k-2}}$. On the other hand, when $X_{k-1}^{\text{loop}} = X_{\zeta_{k-2}+1}$ we know that $X_{\zeta_{k-2}}$ has moved on $\mathcal{P}$ implying that $\zeta_{k-1} = \zeta_{k-2}+1$. Consequently, regardless the finiteness of $\zeta_{k-1}$ -- since $W_k$ is independent of the \textit{BGRW}, $\mathcal{B}_k$ is obtained by adding a loop on $X^{\text{loop}}_{k-1}$ independently of the whole past and with probability $p$. 
	
	To define $X_k^{\text{loop}}$ we proceed in the following way
	\begin{itemize}
		\item If $\zeta_{k-1} < \infty$, 
		\begin{itemize}
			\item  we set $X_k^{\text{loop}} = X_{\zeta_{k-1}+1}$, if $X_{\zeta_{k-1}}$ moves on $\mathcal{P}$;
			\item otherwise, when $X_{\zeta_{k-1}}$ jumps outside $\mathcal{P}$, we let $X_k^{\text{loop}}$ be $X_{\zeta_{k-1}}$. 
		\end{itemize}
		\item In case $\zeta_{k-1} = \infty$, we select uniformly an edge of $X_{k-1}^{\text{loop}}$ on $\mathcal{B}_k$ to walk $X_{k-1}^{\text{loop}}$ through.
	\end{itemize}
	
	By definition of the above coupling, we obtain that $\eta_{y} \ge \eta_0^{\text{loop}}$. The best scenario would be that in which the \textit{BGRW} only walks over $\mathcal{P}$, in this case the stopping times are equals. This concludes the proof.
\end{proof}
As a straightforward consequence of the above coupling, we restate Corollary \ref{cor:xlooptime} in terms of the hitting time $\eta_{y}$.
\begin{corollary}\label{cor:etaybound} Let $T_0$ be a rooted locally finite tree, $x_0$ one of its vertices different from the root and $y$ an ancestor of $x_0$ at distance $\ell$. Then there exists a constant $C$ depending on $p$ only, such that 
\begin{equation*}
\Pd_{T_0,x_0,p}\left( \eta_{y} \le e^{\sqrt{\ell}}\right) \le \frac{C}{\sqrt{\ell}}.
\end{equation*}
\end{corollary}

\section{Positive drift and its consequences}\label{sec:positivedrift}

In this section we show that the \textit{BGRW} has a positive drift away from the root. In particular, we show that 
\[
	\liminf_{n \rightarrow \infty}\frac{\mathrm{dist}_{T_n}(X_n,\text{root})}{n} >0,
\]
almost surely, which implies the \textit{transience} of the walker $X$. We do that by tracking the distance of $X$ from the \textit{root} at random times and comparing it to a right biased simple random walk on $\Z$. Furthermore, this comparison with the right biased simple random walk allows us to improve the results given in Proposition~\ref{prop:harderpart} and Corollary~\ref{cor:etaybound}. Specifically, we can now prove
\begin{itemize}
	\item the walker achieves distances of order $n$ in $n$ steps \textit{w.h.p};
	\item the probability of the walker going back long distances is exponentially rare.
\end{itemize}
Intuitively, the main message of this section is that if we look at the distance of $X$ from the root properly normalized and at some ``random times" we see a random walk on the line that dominates a right biased simple random walk. Let us begin by formally define what we mean by ``random times". For a fixed positive integer $r$, we define three stopping conditions from a vertex $x_0$.
\begin{enumerate}
	\item $\mathrm{dist}_{T_n}(X_n, \text{root}) - \mathrm{dist}_{T_0}(x_0, \text{root}) = r$;
	\item $\mathrm{dist}_{T_n}(X_n, \text{root}) - \mathrm{dist}_{T_0}(x_0, \text{root}) = -r$;
	\item $X$ walks $\exp\{{\sqrt{r}}\}$ steps and none of the previous conditions occurs.
\end{enumerate}
We say $(1)-(3)$ occurred from $X_m$ whenever one of the three stopping conditions occurred when we put $x_0$ as $X_m$. We must point out that if $\mathrm{dist}_{T_m}(X_m, \text{root}) < r$, then stopping condition~$(2)$ cannot be attained.

We define our sequence of stopping times as follows:
\begin{equation}
\begin{split}
	\sigma_0 \equiv & 0; \\
	\sigma_k := & \inf \left \lbrace m > \sigma^{(r)}_{k-1} \mid (1) - (3) \text{ occurs from }X_{\sigma_{k-1}}\right \rbrace.
\end{split}
\end{equation}
From the definition of $(1)-(3)$ follows that, for all $k$, $\sigma^{(r)}_k$ is bounded from above by $\exp\{{\sqrt{r}}\}k$. To avoid clutter, when $r$ is fixed we suppress the upper script from the nation of $\sigma^{(r)}_k$.
\begin{lemma}[Coupling to the biased random walk]\label{lema:couplingsk} Let $T_0$ be a rooted locally finite tree, $x_0$ one of its vertices and $\left \lbrace S_k \right \rbrace_{k \ge 0}$  a right biased simple random walk on $\Z$. Then, there exists a large enough $r = r(p)$ depending on $p$ only, such that the process $\{ \mathrm{dist}_{T_{\sigma_k}}(X_{\sigma_k}, \text{root})/r\}_{k \ge 0}$ starting from $(T_0,x_0)$ and $\left \lbrace S_k \right \rbrace_{k \ge 0}$ starting from $\lfloor \mathrm{dist}_{T_{0}}(x_0, \text{root})/r \rfloor$ can be coupled in such way that
\[
\Pd \left( \mathrm{dist}_{T_{\sigma_k}}(X_{\sigma_k}, \text{root}) \ge rS_k, \forall k \right) = 1.
\]
\end{lemma}
\begin{proof} We begin by a few notations. To simplify our writing, put $d_k := \mathrm{dist}_{T_{\sigma_k}}(X_{\sigma_k}, \text{root})$. Note that the process $\left \lbrace d_k \right \rbrace_{k \ge 0}$ is a random walk on the half line whose increments belong to the interval $[-r,r]$. Let $\left \lbrace U_k \right \rbrace_{k \ge 0}$ be a sequence of \textit{i.i.d} random variables independent of the \textit{BGRW} and following uniform distribution on~$[0,1]$. Also let $\tilde{\mathcal{F}}_k$ be the $\sigma$-algebra generated by the \textit{BGRW} process up to time $\sigma_k$ and let~$\mathcal{H}_k$ be the $\sigma$-algebra encoding all information of $\tilde{\mathcal{F}}_k$ and all the uniform random variables up to time $k$.

Regarding the increments of $\left \lbrace d_k \right \rbrace_{k \ge 0}$, we claim
\begin{claim}\label{claim:q} There exist positive constants $r$ and $C$ depending on $p$ only, such that for all $k$ 
	\begin{equation}
	\inf_{T_0,x_0} \Pd_{T_0,x_0,p} \left( \Delta d_{k+1} = r \middle | \tilde{\mathcal{F}}_k\right) > 1-\frac{C}{\sqrt{r}}.
	\end{equation}
\end{claim}
\begin{claimproof}{} First of all, observe that $\Delta d_{k+1} = r$ if, and only if, $\sigma_{k+1}$ stops because of condition $(1)$. Also, by the strong Markov Property, it is enough to prove the claim for $\sigma_1$. Said that, let us first assume we start from an initial condition $(T_0,x_0)$ with~$\dist{0}{x_0}{\text{root}} \ge r$. We derive the desired lower bound by proving upper bounds for the probability of the events~$\{\Delta d_{1} = -r\}$ and~$\{|\Delta d_{1}|<r \}$. The former occurs if, and only if, $\sigma_1$ stops because condition $(2)$ and the latter because of condition $(3)$. 
	
Observe that if $\sigma_1$ stops because of $(2)$ then the walker $X$ visited the ancestor $y$ of $x_0$ at distance $r$, since $\dist{0}{x_0}{ \text{root}}\ge r$, spending at most $\exp\{\sqrt{r}\}$ steps. So, using Corollary~\ref{cor:etaybound} we have
\begin{equation}
\Pd_{T_0,x_0,p}\left( \sigma_1 \textit{ stops because of } (2) \right) \le \Pd_{T_0,x_0,p}\left( \eta_y \le e^{\sqrt{r}} \right) \le \frac{C'}{\sqrt{r}}.
\end{equation} 
Note that the above upper bound holds for all possible pairs of a rooted locally finite tree $T_0$ and one of its vertex $x_0$ at distance greater than $r$ from the root.
 
Finally, if $X$ stops because of the occurrence of $(3)$ then $X$ has walked for $e^{\sqrt{r}}$ steps and has not visited the ancestor $y$ of $x_0$ at distance $r$ neither has increased its distance from $y$ by $r$. This is the same that we observe a process $\tilde{X}$ on the subtree $T_y$ that in $e^{\sqrt{r}}$ steps does not be at distance $2r$ from the root $y$. Applying Proposition \ref{prop:harderpart} (page \pageref{prop:harderpart}) with $n = e^{\sqrt{r}}$ and~$M=4$ we obtain that there exist $r_0 = r(p)$ such that for all $r \ge r_0$ we have
\begin{equation}
\inf_{T_0,x_0} \Pd_{T_0,x_0,p} \left( \exists m \le e^{\sqrt{r}}, \dist{m}{\tilde{X}_m} {\tilde{\textit{root}}} > r^2\right) \ge 1-  \exp\{-e^{\sqrt{r}/4}\}.
\end{equation}
Choosing $r$ large enough so that $r \ge r_0$, $\exp\{-e^{\sqrt{r}/4}\} \le C'/\sqrt{r}$ and  $r^2 > 2r$ we conclude that
\begin{equation}
\sup_{T_0,x_0}\Pd_{T_0,x_0,p}\left( \sigma_1 \textit{ stops because of } (3) \right) \le \exp\{-e^{\sqrt{r}/4}\}.
\end{equation}
To drop our assumption that we start at distance greater than $r$ from the root, just recall that when this is not the case the condition $(2)$ has probability zero and the above upper bound for condition $(3)$ still holds. This implies that
\begin{equation}
\inf_{T_0,x_0} \Pd_{T_0,x_0,p}\left( \sigma_1 \textit{ stops because of } (1) \right) \ge 1 - \frac{C}{\sqrt{r}}
\end{equation}
which combined with strong Markov Property proves the claim.
\end{claimproof}

Since we may increase $r$, we choose it large enough so 
\begin{equation}
q := 1-\frac{C}{\sqrt{r}} > \frac{1}{2}.
\end{equation}

Now we couple the process $\{d_k/r\}_{k\ge1}$ and $\{S_k\}_{k\ge1}$ in the following way. Set
\begin{equation}
S_0 := \left \lfloor \frac{\mathrm{dist}_{T_0}(x_0, \text{root})}{r} \right \rfloor
\end{equation}
and assume we have defined $\{S_j\}_{j=0}^{k-1}$ in such way that it has the distribution of $k-1$ steps of the desired right biased random walk. Let $Q_{k-1}$ denote
\begin{equation}\label{def:qk}
Q_{k-1} := \Pd_{T_0,x_0,p}( \Delta d_{k} = r| \tilde{\mathcal{F}}_{k-1} )
\end{equation}
and recall from Claim \ref{claim:q} that $Q_{k} \ge q$ for all $k$. We then define $S_k$ in the following way
\begin{equation}
S_k = \left\{
\begin{array}{ll}
S_{k-1}+1  & \mbox{if } \Delta d_k = r \text{ and } U_k \le q/Q_{k-1},  \\
S_{k-1} - 1 & \mbox{otherwise. }
\end{array}
\right.
\end{equation}
In words, if at time $\sigma_k$ the walker $X$ increased its distance from the root by $r$, then $S_k$ jumps to the right with probability  $q/Q_{k-1}$. In this way, whenever the process $\{d_k/r\}_{k\ge 0}$ jumps back (at most one unit), the \textit{SRW} also jumps back one unit.

Now, we show that the process $\{S_k\}_{k\ge 0}$ does have the distribution we desire. We start by checking that the increments are $1$ with probability $q$ or $-1$ with probability $1-q$. By the definition of $S$ we have
\begin{equation}\label{eq:probsk}
\begin{split}
\Pd_{T_0,x_0,p}\left( \Delta S_{k} = 1 \middle | \mathcal{H}_{k-1}\right)  & = \Pd_{T_0,x_0,p}\left( \Delta d_{k+1} = r, U_k \le q/Q_{k-1} \middle | \mathcal{H}_{k-1}\right) \\
& = \Pd_{T_0,x_0,p}( \Delta d_{k} = r| \mathcal{H}_{k-1} )\frac{q}{Q_{k-1}} \\
& = q,
\end{split}
\end{equation}
since $U_k$ is independent of $\mathcal{H}_{k-1}$ and~$Q_{k-1}$ is measurable with respect to $\mathcal{H}_{k-1}$. Moreover, the equality below holds
\[
Q_{k-1} = \Pd_{T_0,x_0,p}( \Delta d_{k} = r| \tilde{\mathcal{F}}_{k-1} ) = \Pd_{T_0,x_0,p}( \Delta d_{k} = r| \mathcal{H}_{k-1} )
\]
since $\mathcal{H}_{k-1}$ is $\tilde{\mathcal{F}}_{k-1}$ added of information independent of the whole process $\{(T_k,X_k)\}_{k\ge1}$. Equation (\ref{eq:probsk}) allows us to derive the independence of all the increments of our right biased random walk $\{S_k\}_{k \ge 1}$. For any fixed set of indexes $k_1 < k_2 < \cdots < k_j$ and any vector~$(a_1,\cdots,a_j) \in \{-1,1\}^j$ we have
\begin{equation*}
\begin{split}
\Pd\left( \Delta S_{k_1} = a_1, \cdots, \Delta S_{k_j} = a_j \right) & = \Ed \left [ \mathbb{1}_{\{ \Delta S_{k_1} = a_1 \}} \cdots \mathbb{1}_{\{ \Delta S_{k_{j-1}}=a_{j-1} \}} \Ed \left[ \mathbb{1}_{\{ \Delta S_{k_j} = a_j \}}\middle | \mathcal{H}_{k_j-1}\right] \right] \\
& = \Pd\left( \Delta S_{k_1} = a_1, \cdots, \Delta S_{k_{j-1}} = a_{j-1} \right)\Pd\left(\Delta S_{k_j} = a_j \right)
\end{split}
\end{equation*}
which implies the independence of the increments $\Delta S_k$. So, $\{S_k\}_{k\ge 0}$ is distributed as simple random walk on $\Z$ with probability $q > 1/2$ of jumping to the right starting at the left of~$d_0/r$. And, by construction, we have that $d_k/r \ge S_k$ for all $k$. This concludes the proof of the lemma.
\end{proof}
As a consequence of the coupling above, we prove the transience of the walker. 
\begin{proposition}[Transience of the walker]\label{prop:liminf} There exists a constant $c >0$ depending on $p$ only, such that
	\[
	\liminf_{n \rightarrow \infty}\frac{\mathrm{dist}_{T_n}(X_n,\mathrm{root})}{n}  \geq c, \; \Pd_{T_0,x_0,p}\textit{- almost surely},
	\]
for all initial conditions $(T_0,x_0)$.
\end{proposition}
\begin{proof} Lemma \ref{lema:couplingsk} guarantees the existence of a positive constant $r$, depending on $p$ only, so that, for any initial condition on the \textit{BGRW}, 
	\begin{equation}
	\mathrm{dist}_{T_{\s_k}}(X_{\sigma_k}, \text{root}) \ge r S_k, \textit{a.s.,}
	\end{equation}
	where $\{S_k\}_{k\ge 0}$ is a right biased simple random walk on $\Z$. By the Strong Law of Large Numbers we also have that
	\begin{equation}
	\lim_{k \rightarrow \infty} \frac{S_k}{k} = \mu(r) = \mu\textit{, a.s..}
	\end{equation}
	This implies that
	\begin{equation}\label{eq:liminfx}
	\liminf_{k \rightarrow \infty} \frac{\mathrm{dist}_{T_{\s_k}}(X_{\sigma_k}, \text{root}) }{k} \ge r\mu \textit{, a.s..}
	\end{equation}
	On the other hand, by the definition of the stopping times $\sigma_k$ we have that, for all $k$, the following inequality holds
	\[
	|\sigma_{k+1} - \sigma_k| \le \exp\{\sqrt{r}\},
	\]
	almost surely, as well as
	\[
	\left |\mathrm{dist}_{T_{\s_k}}\left(X_{\sigma_k}, \text{root}\right) - \mathrm{dist}_{T_{\s_{k-1}}}\left(X_{\sigma_{k-1}}, \text{root}\right) \right | \le r,
	\]
	almost surely. Thus, if $n$ is an integer such that $ n \in [\sigma_k, \sigma_{k+1}]$ we have
	\begin{equation}
	\frac{\mathrm{dist}_{T_n}\left(X_n, \text{root}\right)}{n} \ge \frac{\mathrm{dist}_{T_{\s_{k+1}}}\left(X_{\sigma_{k+1}}, \text{root}\right) - r}{e^{\sqrt{r}}(k+1)} 
	\end{equation}\label{ineq:posliminf}
	which combined with (\ref{eq:liminfx}) proves the that
	\begin{equation}
	\liminf_{n \rightarrow \infty}\frac{\mathrm{dist}_{T_n}(X_n,\mathrm{root})}{n}  = c >0, \textit{a.s},
	\end{equation}
	for some positive constant $c$ depending on $p$ only.
\end{proof}
\subsection{Controlling returns and degrees}
The coupling gives us a picture of the dynamics of the walker: it is moving away from the root at linear speed. Before moving on, we use the coupling to show that vertices are not visited many times and (as a result) degrees in our tree tend to be small.

 \begin{lemma}\label{lemma:returntime} Let $(T_0,x_0)$ be the initial state of the dynamics and $y\in T_0$ be given. Then for any $0<p\leq 1$, the number $v(y)$ of visits to $y$,
\[v(y):=|\{0\leq t<+\infty\,:\, X_t=y\}|\]
satisfies \[\forall k\geq 0\,:\, \Pd_{T_0,x_0,p}(v(y)\geq k)\leq C\,e^{-\beta\,k}\]
and 
\[\Pd_{T_0,x_0,p}(v(y)>0)\leq C\,e^{-\beta\,{\rm dist}_{T_0}(x_0,y)}\] 
for some $\beta,C>0$ depend only on $p$. Moreover, if $p_0>0$ and $p_0\leq p\leq 1$, $C$ and $\alpha$ can be chosen uniformly in $p_0$. \end{lemma}
\begin{proof}
	Recall from the coupling we just constructed that we may choose $r$ large enough so that the process $\{\dist{\s_k}{X_{\s_k}}{y}/r\}_k$ dominates a right biased random walk $\{S_k\}_k$ on $\Z$ starting on $s_0 := \lfloor\frac{\dist{0}{x_0}{y}}{r} \rfloor$. We count the visits to $y$ per time interval $(\s_{k-1},\s_k]$ (with a possible additional visit at time $0$).
	\[v(y):= \Ind{\{X_0=y\}}  + \sum_{k=1}^{\infty}\sum_{t=\s_{k-1}+1}^{\s_k}\Ind{\{X_t=y\}}.\]
	Recalling that $\s_k-\s_{k-1}\leq e^{\sqrt{r}}$, we may rewrite the expression as:
	\[v(y)\leq \Ind{X_0=y}  + e^{\sqrt{r}}\sum_{k=1}^{\infty}\Ind{\{X_t=y\text{ for some }\s_{k-1}<t\leq \s_k\}}.\]
	Now, by the definition of the stopping time sequence $\{\s_{k+1}\}_k$, for $t\in [\s_k,\s_{k+1}]$ then~\[|\dist{\s_k}{X_{\s_k}}{y} - \dist{t}{X_{t}}{y}| \le r.\]
	As a result, if $X_t=y$ for some $\s_{k-1}\leq t\leq \s_k$, then $S_k\leq \dist{\s_k}{X_{\s_k}}{y}/r\leq 1$. In particular, 
	\[v(y)\leq \Ind{\{X_0=y\}}  + e^{\sqrt{r}}\sum_{k=1}^{\infty}\Ind{\{S_k\leq 1\}}.\]
The RHS is (up to a constant) the number of visits of a right-biased random walk started from $S_0\geq 0$ to the interval $(-\infty,1]$. This number has an exponential tail that only depends on the bias. The first inequality follows because we can choose $r=r(p_0)$ to guarantee a bias of $1/3$ (say) for all $p_0\leq p\leq 1$. The second inequality also follows (perhaps with changes to $C,\beta$) once we realize that $S_0\geq \dist{0}{x_0}{y}/r$.\end{proof}

\begin{lemma}\label{lem:degreetail}Given $p_0>0$, there exist constants $C,\alpha>0$ depending only on $p_0$ such that, for all $p_0\leq p\leq 1$, all finite $(T_0,x_0)$, and all $n,k\geq 1$,
\[
\sum_{t=0}^n\Pd_{T_0,x_0,p}(\degr{t}{X_t}\geq k)\leq C\,(n+|V(T_0)|)\,e^{\alpha\,(\Delta(T_0)-k)}.
\]
\end{lemma}
\begin{proof}Assume without loss of generality that $V(T)=\{1,\dots,\ell\}$. Also let ${\ell+1},{\ell+2},\dots$ be the vertices that the BGRW process started from $T,x$ creates. Finally, we let $v_t(i)$ denote the number of visits to $i$ up to time $t$, so that $v_t(i)=0$ if $i>t+\ell$. Note that $d_t(v_i)\geq k$ implies that $v_n(i)\geq k-\Delta(T)$, as the degree only grows at times when $i$ is visited. Therefore, for all $0\leq t\leq n$, $\alpha>0$ and $k\geq \ell+1$:
\begin{equation*}
\begin{split}
\Pd_{T,o,p}(\degr{t}{X_t}\geq k) = & \sum_{i=1}^{\ell+n}\Pd_{T,o,p}(X_t=i,\degr{t}{i}\geq k)\\
&\leq e^{\alpha\,(\Delta(T)-k)}\,\sum_{i=1}^{\ell+n}\Ed_{T,o,p}[\Ind{\{X_t=i\}}\,\exp(\alpha\,v_n(i))].
\end{split}
\end{equation*}
As a consequence:
\begin{eqnarray*}
	\sum_{t=0}^n\Pd_{T,o,p}(\degr{t}{X_t}\geq k) &\leq & e^{\alpha\,(\Delta(T)-k)}\,\sum_{t=0}^n\sum_{i=1}^{\ell+n}\Ed_{T,o,p}[\Ind{\{X_t=i\}}\,\exp(\alpha\,v_n(i))]\\ \mbox{(use $\sum_{t\leq n}\Ind{\{X_t=i\}}=V_n(i)$)} &\leq & e^{\alpha\,(\Delta(T)-k)}\sum_{i=1}^{\ell+n}\Ed_{T,o,p}[v_n(i)\,\exp(\alpha\,v_n(i))].\end{eqnarray*}
By Lemma \ref{lemma:returntime}, $v_n(i)\leq v(i)$ has an exponential tail uniformly in $i,n,T_0,x_0$ and $p\in [p_0,1]$. Thus we may take $\alpha=\beta/2$ (with $\beta$ from that Lemma) and adjust $C$ to obtain our result.\end{proof}

\section{The local point of view on infinite trees}\label{sec:finiteinfinite}


In the last section we obtained that the distance between the random walker and its starting position grows at least like $n$. Obtaining sharper results will require a deeper understanding of the process. This section makes some progress in this direction by introducing the right state space for this purpose.

\subsection{Rooted graphs and trees}

We recall the definition of rooted graphs, rooted trees and the local topology on these objects. All notions we need are defined and studied in Bordenave's lecture notes \cite{Bordenave2016}.

A {\em rooted graph} is a pair $(G,o)$ where $G$ is a connected, locally finite graph and $o$ is vertex of $G$. Note that any rooted graph must have a countable vertex set which we may assume to be a subset of $\N$ or $\Z$. Two rooted graphs $(G,o)$ and $(G',o')$ are rooted-isomorphic -- denoted by $(G,o)\simeq (G,o')$ -- if there is a bijection of the vertex sets of $G$ and $G'$ that maps $o$ to $o'$ and preserves edges. We let $[G,o]$ denote the {\em equivalence class} of $(G,o)$ under rooted isomorphisms, and let $\sG_*$ denote the set of all such $[G,o]$. 

Given $r\in\N$ and a rooted graph $(G,o)$, $(G,o)_r$ denotes the graph obtained by retaining only the vertices of $G$ within graph-theoretic distance $r$ from $o$ and the edges between those vertices. Clearly, $(G,o)\simeq (G',o')$ implies $(G,o)_r\simeq (G'o')_r$  for all $r>0$, and we may define $[G,o]_r$ as the equivalence class of $(G,o)_r$. The {\em distance} between $[G,o],[G',o']\in\sG_*$ is defined by:
\[\rho([G,o],[G',o']):=\frac{1}{1+\sup\{r\in\N\,:\, [G,o]_r = [G',o']_r\}}.\]
One can show that $(\sG_*,\rho)$ is a Polish metric space. 

The set $\mathcal{T}_*\subset\sG_*$ of rooted trees is the set of equivalence classes $[T,o]$ where $T$ is a locally finite tree. This is a closed subset of $\sG_*$ and is therefore a Polish space with the metric $\rho$.


Our \textit{BGRW} dynamics (see, Section~\ref{sec:markov-kernel}) may be naturally extended to the set $\sT_*$ of random rooted trees. The idea is that the state $[T_t,X_t]$ describes the tree ``rooted at the position of the walker". With some abuse of notation, we use $K_p$ to denote the Markov transition kernel of our process over this space as well.


\subsection{Empirical measures and local functions}\label{sec:localweak}

Much of the remainder of the paper will be spent dealing with the tree as viewed by the walker. More precisely, we will study the {\em empirical measure of the tree around the walker}.

\begin{definition}\label{def:empiricalmeasure} Given a realization $[T_t,X_t]$ of the process $K_p$, we let $\widehat{P}_n$ denote the {\em empirical measure}, that is, the random probability measure over $\sT_*$ given by:
\[\widehat{P}_n = \frac{1}{n+1}\sum_{t=0}^{n}\delta_{[T_t,X_t]}.\]\end{definition}

Thus for a given element $[T,v]\in\sT_*$ and $r\in\N$,
\[(n+1)\widehat{P}_n([G,o]\in\sT_*\,:\, [G,o]_r = [T,v]_r)\]
counts the number of times $0\leq t\leq n$ at which the ball of radius $r$ around $X_t$ in $T_t$ is isomorphic to $[T,v]_r$. 

In this section we show that, under suitable assumptions on $[T_0,x_0]$, the walker has a well-defined positive speed and in Section~\ref{sec:localweak} that   $\widehat{P}_n\Rightarrow P_p$ almost surely, where $P_p$ is an invariant measure for the process $K_p$. In the proof of such results, {\em local functions} will play an important role.

\begin{definition}[Local function] A function $\psi:\sT_*\to\R$ is said to be $r$-local if $\psi([T,o]) = \psi([T,o]_r)$ for all $[T,o]\in \sT_*$. A function $\psi$ is local if it is $r$-local for some $r\in\N$. 
\end{definition}
To avoid clutter, with a slight abuse of notation, we will write $\psi(G,o)$ instead of $\psi([G,o])$.

The main result of this section is about the convergence of empirical averages of \textit{bounded local functions} on the space~$\sT_*$ along trajectories. In order to do that we restrict ourselves on measures on $\sT_*$ which have some nice properties. To define this class of measures we first define certain stopping times. Given $\ell\in\N\backslash\{0\}$, let $Q_\ell$ be the path of length $\ell$, i.e. the graph with vertex set $\{0,1,\dots,\ell\}$ and edges $\{(i-1),i\}$ ($1\leq i\leq \ell$). We define $\tau_\ell$ as the first time when the graph around $X_t$ in $T_t$ is $Q_\ell$.
\[\tau_\ell :=\inf\{t\geq 0\,:\, [T_t,X_t]_\ell = [Q_\ell,\ell]\}.\]

\begin{definition}A measure $\mu$ over pairs $(T,x)$ is called $p$-escapable if $\tau_{\ell}<+\infty$ $\Pd_{\mu,p}$-almost surely for every $\ell\geq 1$.\end{definition}
We can now state the main theorem of this section.
\begin{theorem}\label{thm:rfuncthm}For each parameter $0<p\leq 1$ and each bounded integer $r \ge 1$ and all $r$-local function $\psi:\sT_*\to \R$, there exists a constant $\mathcal{M}_p\psi$ such that
	\[
	\lim_{n\to +\infty}\widehat{P}_n\psi =\mathcal{M}_p\psi\mbox{ }\mathbb{P}_{\mu,p }\mbox{-a.s.},
	\]
	for all $p$-escapable distribution $\mu$ on $\sT_*$.
\end{theorem}
As we will see in Section~\ref{sec:localweak}, this result implies that ``the tree seen by the walker"~converges to an infinite random tree when $n\to +\infty$.  

The proof of the theorem relies on the following results.
The first Lemma shows that the empirical measures of processes started from two $p$-escapable measures are asymptotically the same, i.e. the initial distribution is ``forgotten". The second Lemma shows that the empirical measures of local functions when  our process starts from specific   $p$-escapable measures, corresponding to $\delta_{[T_0, x_0]}$ with $T_0$ finite, converges to a constant. 
 
\begin{lemma}[Forgetfulness of empirical measures; proof in \S \ref{ss:forget}]\label{lem:forget} Let $\mu$, $\nu$ be two $p$-escapable measures. Let $r$ be a positive integer and $\psi:\sT_*\to\R$ be a bounded $r$-local function. If there exists a constant~$\mathcal{M}_p\psi$ such that $\widehat{P}_n\psi\to \mathcal{M}_p\psi$ $\Pd_{\nu,p}$-almost-surely, then we also have $\widehat{P}_n\psi\to \mathcal{M}_p\psi$ $\Pd_{\mu,p}$-almost-surely.
\end{lemma}
The above statement allows us to restrict ourselves to subclasses of initial distributions and then extend the results to the whole class of $p$-escapable distributions. This procedure greatly reduces the amount of work in our proofs. The next lemma is a finite version of Theorem~\ref{thm:rfuncthm}.
\begin{lemma}[proof in \S \ref{ss:rfunclema}]\label{lema:rfunclema}Consider an initial condition $(T_0,x_0)$ with $T_0$ finite. Fix $0<p\leq 1$ and a bounded $r$-local function $\psi:\sT_*\to \R$ (with $r$ a positive integer). Then there exists a constant $\mathcal{M}_p\psi$ such that $\widehat{P}_n\psi \to \mathcal{M}_p\psi, \mbox{   } \mathbb{P}_{T_0,x_0,p }$\mbox{-a.s.}
\end{lemma}
Combining the two lemmas gives us a straightforward proof of Theorem~\ref{thm:rfuncthm}.
\begin{proof}[Proof of Theorem~\ref{thm:rfuncthm}]
Let $[T_0,x_0]$ be a finite rooted tree and consider $\nu = \delta_{[T_0, x_0]}$. By Lemma~\ref{lema:rfunclema}, we have that~$\widehat{P}_n\psi \to \mathcal{M}_p\psi, \mbox{   } \mathbb{P}_{\nu,p }$\mbox{-a.s.} But, by Lemma~\ref{lem:degreetail}, $\nu$ is $p$-escapable. Therefore, by Lemma~\ref{lem:forget},  $\widehat{P}_n\psi \to \mathcal{M}_p\psi, \mbox{   } \mathbb{P}_{\mu,p }$\mbox{-a.s.}
\end{proof}
Now we know how the main theorem follows from our lemmas we organize the remainder of this section as follows. In the next subsection we prove Theorem~\ref{thm:positivespeed} as an application of Theorem~\ref{thm:rfuncthm}. Then, in Subsection~\ref{sec:escapable}, we discuss the concept of $p$-escapable trees and give quantitative criteria for escapabality. Finally, in Subsections~\ref{ss:forget} and \ref{ss:rfunclema}, we prove lemmas~\ref{lem:forget} and \ref{lema:rfunclema}, respectively.

\subsection{Existence of the speed} In this section we prove a stronger version of Theorem~\ref{thm:positivespeed} as an application of Theorem~\ref{thm:rfuncthm}. We prove that the random walker in the~\textit{BGRW} model has a well-defined positive speed for any $p$-escapable initial distribution on~$\sT_*$.
\begin{theorem}[Linear speed of the walker] For each $0<p\le 1$, there exists a constant $0<c(p)\le p$ such that, for any $p$-escapable distribution $\mu$ on $\sT_*$,
	\[
	\lim_{n\rightarrow \infty } \frac{\dist{n}{X_n}{\mathrm{root}}}{n}  = c(p) \; \Pd_{\mu,p}\textit{-a.s.}.
	\]
	
\end{theorem}
\begin{proof}Our strategy is to prove that the distance from the root at time $n$ may be written as a sum of a bounded term, a martingale, and a sum of a $1$-local function computed along the trajectory. The result will then follow from Theorem~\ref{thm:rfuncthm}. 

Fix a $p$-escapable measure $\mu$ and define the following bounded $1$-local function:
\[
	\psi(T,x) :=  \frac{p(\degr{}{x}-1)}{\degr{}{x} + 1} + \frac{(1-p)(\degr{}{x} - 2)}{\degr{}{x}}.\]

 Note that:
	\begin{equation}
	\begin{split}
	\mathbb{E}_{\mu,p}\left[ \Delta \dist{n}{X_n}{ \mathrm{root}} \mid  \mathcal{F}_{n}\right]  
	& = \psi(T_n,X_n) + (1 - \psi(T_n,X_n))\,\Ind{\{X_n=\mathrm{root}\}}.
	\end{split}
	\end{equation}
	Thus 
	\begin{equation}
	\dist{n+1}{X_{n+1}}{\mathrm{root}} = M_{n+1} + \sum_{j=0}^n \psi(T_j,X_j ) + C_{n+1},
	\end{equation}
	where $\{M_j\}_j$ is a martingale with bounded increments and
	\[C_{n+1}:= \sum_{j=0}^n(1 - \psi(T_n,X_n))\,\Ind{\{X_n=\mathrm{root}\}}\]
	is bounded by the total number of visits to the root, which is a.s. finite because our process is transient. 
	
	Since the martingale has bounded increments, Azuma's inequality implies that $M_n/n\to 0$ almost surely. Since $C_{n}$ is a.s. bounded, $C_{n}/n\to 0$ almost surely as well. Theorem \ref{thm:rfuncthm} gives us that 
	\begin{equation}\label{eq:convpsi}
	\lim_{n\to +\infty}\frac{\dist{n+1}{X_{n+1}}{\mathrm{root}}}{n+1} = \lim_{n\rightarrow \infty }\frac{1}{n+1}\sum_{j=0}^n \psi(T_j,X_j) = c(p):=\mathcal{M}_p(\psi), \; \Pd_{\mu,p}\textit{-a.s},
	\end{equation}
	with $\mathcal{M}_p(\psi)$ is a constant depending on $p$ and $\psi$ only (but not on $\mu$). Proposition~\ref{prop:liminf} implies $c(p)>0$. \end{proof}

\subsection{Escapable trees and mixing}\label{sec:escapable}

The concept of escapable trees is closely related to the forgetfulness of the walker. Roughly speaking, once the walker has escaped it may never come back to its initial tree. In this section we prove a quantitative criteria for checking ``escapability".  Let us notice that not all trees are $p$-escapable. The next example gives a useful example to keep in mind.

\begin{example}[A hard tree]\label{ex:hardtree}Consider an infinite rooted tree $T_{\rm hard}$, with a root node $x_{\rm hard}$, such that each node at distance $h$ from the root has $g(h)>0$ children, with $1/g(h)$ summable. Consider the BGRW from $(T_{\rm hard},x_{\rm hard})$. One can show via the Borel-Cantelli Lemma that there is a positive probability that ${\rm dist}_{T_t}(X_t,x_{\rm hard}) = t$ for all $t\in\N$ (i.e. the walker simply ``walks down the tree"~for eternity). For the same reason, for any $\ell\in\N\backslash\{0\}$ there is a positive probability that $\tau_{\ell}=+\infty$. \end{example}

The next Lemma (proven subsequently) explains one reason why $p$-escapability is important to us. The  Lemma gives a quantitative criterion for ``escapability". In a way, it says that what makes the tree in Example \ref{ex:hardtree} non-escapable is that $X_t$ sees very large degrees along the way. 

\begin{lemma}[Quantitative escapability]\label{lem:quantitative} Assume $\mu$ is a starting measure for which there exist $C,D,\alpha>0$ such that, for all $n,k\geq 1$,
	\begin{equation}\label{eq:condicaolegal1}\sum_{t=0}^n\Pd_{\mu,p}(d_{T_t}(X_t)\geq k)\leq C\,(n+D)\,e^{-\alpha\,k}.\end{equation}
	Then $\mu$ is $p$-escapable and
	\begin{equation}\label{eq:estimateescapemu}\Ed_{\mu,p}[\tau_\ell\wedge n]\leq \frac{\ell 2^{\ell-1}}{p^\ell}\left\lceil \frac{2\log(n+D) + \log C}{\alpha}\right\rceil + 1.\end{equation}
	Moreover, the following event holds with probability $\geq 1 - 1/n^{16}$: for all times $1\leq t\leq n$,  
	\begin{equation}\label{eq:estimateescapefinite}\Ed_{T_t,X_t,p}[\tau_\ell\wedge n]\leq \frac{\ell 2^{\ell-1}}{p^\ell} \left\lceil \frac{\log(2n+D) +17\log n + \log C}{\alpha}\right\rceil + 1.\end{equation}
	In particular, all measures $\mu=\delta_{(T,x)}$, with $T$ a finite tree, satisfy the above inequalities, with~$C\leq e^{\alpha\,\Delta(T)}$ by virtue of Lemma \ref{lem:degreetail}.\end{lemma}
So Lemma \ref{lem:quantitative} also guarantees that not only $\mu$ is $p$-escapable, but also that, with high probability, each intermediate $[T_t,X_t]$ satisfies a quantitative escapability property. One important point is that the assumptions of Lemma \ref{lem:quantitative} always apply to initial trees that are finite, by virtue of Lemma \ref{lem:degreetail}.

\begin{proof} We first prove inequality (\ref{eq:estimateescapemu}) for $\Ed_{\mu,p}[\tau_\ell\wedge n]$, noting that it implies $p$-escapability because, assuming it, we can show:
	\[\Pd_{\mu,p}(\tau_{\ell}=+\infty)\leq \frac{\Ed_{\mu,p}[\tau_\ell\wedge n]}{n}\leq \frac{1}{n}\left(\frac{\ell 2^{\ell-1}}{p^\ell}\left\lceil \frac{2\log(n+D) + \log C}{\alpha}\right\rceil + 1\right)\stackrel{n\to +\infty}{\to} 0,\]
	Fix some $\Delta>0$. Let $\eta_\Delta$ be the first time the walker ``sees" a node with degree $\geq \Delta$,
	\begin{equation}
	\eta_\Delta := \inf \left \lbrace m \ge 0 \; \middle | \; \degr{m}{X_m}\geq \Delta\right \rbrace
	\end{equation}
	and consider the event $A_i$ in which at time $i\ell$ the walker creates a new leaf, jumps to it and grows a path of length $\ell-1$ from it in such way that at time $i\ell+\ell-1$ the walker finds itself at the tip of a path of length $\ell$. Now we define another stopping time
	\begin{equation}
	\tilde{\tau_\ell} := \inf \left \lbrace i \ge 0 \; \middle | \; A_i \text{ occurs }\right \rbrace
	\end{equation}
	and note that $\tau_\ell \le \ell\tilde{\tau_\ell}$.
	
	By the definition of the event $A_i$, we have that
	\[
	\Pd_{\mu,p} \left( A_i \middle | \mathcal{F}_{i\ell}\right) = \frac{p}{\degr{i\ell}{X_{i\ell}}}\left(\frac{p}{2}\right)^{\ell-1}.
	\]
	Thus, by the definition of the stopping time $\eta_\Delta$, we also have
	\begin{equation}
	\mathbb{1}\{\eta_\Delta > i\ell \}\Pd_{\mu,p} \left( A_i \middle | \mathcal{F}_{i\ell}\right) \ge \frac{p}{\Delta}\left(\frac{p}{2}\right)^{\ell-1}\mathbb{1}\{\eta_\Delta > i\ell \}.
	\end{equation}
	Using the above inequality and noticing that $\mathbb{1}\{\eta_\Delta > k\ell\} \le \mathbb{1}\{\eta_\Delta > (k-1)\ell\}$, we obtain
	\begin{equation}
	\begin{split}
	\Pd_{\mu,p}\left( \tilde{\tau_\ell} \ge k, \eta_\Delta > k\ell\right) & = \Ed_{\mu,p} \left[ \Ed_{\mu,p} \left[ \mathbb{1}_{A^c_k} \middle | \mathcal{F}_{k\ell}\right]\mathbb{1}\{\eta_\Delta > k\ell \} \prod_{i=0}^{k-1}\mathbb{1}_{A^c_i} \right] \\
	& \le \left(1- \frac{p}{\Delta}\left(\frac{p}{2}\right)^{\ell-1}\right)\Pd_{\mu,p}\left( \tilde{\tau_\ell} \ge k-1, \eta_\Delta > (k-1)\ell\right).
	\end{split}
	\end{equation}
	Then, proceeding by induction, we obtain the following upper bound
	\begin{equation}
	\Pd_{\mu,p}\left( \tilde{\tau_\ell} \ge k, \eta_\Delta > k\ell\right) \le \left(1- \frac{p}{\Delta}\left(\frac{p}{2}\right)^{\ell-1}\right)^k.
	\end{equation}
	Since
	\[
	\Pd_{\mu,p}\left( \tau_\ell \wedge n \ge k \ell \right) \le \left\{\begin{array}{ll}\Pd_{\mu,p}\left( \tilde{\tau_\ell} \ge k \right),
	& k\ell\leq n\\ 0, & \mbox{otherwise};\end{array}\right.
	\]
	we have: 
	\begin{align*}
	\Pd_{\mu,p}\left( \tau_\ell \wedge n \ge k\ell \right) \leq \Pd_{\mu,p}\left(\tilde{\tau_\ell} \ge k ,\eta_\Delta>k\ell\right) + \Pd_{\mu,p}\left(\eta_\Delta\leq k\ell\right) ,
	\end{align*}
	and consequently,
	\begin{equation}\label{eq:expttilde}
	\Ed_{\mu,p}\left[ \tau_\ell\wedge n \right] \le
	\sum_{j=0}^{\lfloor \frac{n}{\ell}\rfloor}
	\ell\Pd_{\mu,p}\left( \tau_\ell \wedge n \ge j \,r \right)
	\le \frac{\ell\,2^{\ell-1}}{p^\ell}\Delta + n\,\Pd_{\mu,p} \left( \eta_\Delta \le n\right).
	\end{equation}
	But, by our condition (\ref{eq:condicaolegal}),
	\begin{equation}\label{eq:conditionused}\Pd_{\mu,p} \left( \eta_\Delta \le n\right)\leq \sum_{t=0}^n\Pd_{\mu,p}(\degr{t}{X_t}\geq \Delta)\leq C\,(n+1)\,e^{-\alpha\,\Delta}.\end{equation}
	so
	\[\Ed_{\mu,p}\left[ \tau_\ell\wedge n \right] \le \frac{\ell 2^{\ell-1}}{p^\ell}\Delta + C\,(n+1)^{2}e^{-\alpha\,\Delta}.\]
	We may choose \[\Delta = \left\lceil \frac{2\log (n+1) + \log C}{\alpha}\right\rceil\] to finish the proof of (\ref{eq:estimateescapemu}). 
	
	To prove the last statement in the Theorem, we go back to (\ref{eq:conditionused}) and note that was the only step in that proof where we used (\ref{eq:condicaolegal1}). In particular, we can go back to (\ref{eq:expttilde}) and obtain:
	\[\Ed_{T_t,X_t,p}\left[ \tau_\ell\wedge n \right] \le \frac{\ell 2^{\ell-1}}{p^\ell}\Delta + n\,\Pd_{T_t,X_t,p}(\eta_\Delta\leq n).\]
	In particular, if $\Delta$ is given and we consider the event  
	\[\mathcal{G}(\Delta):= \bigcap_{t=0}^n\{\Pd_{T_t,X_t,p}(\eta_\Delta\leq n)\leq 1/n)\},\]
	we have that, when $\mathcal{G}(\Delta)$ holds, then:
	\begin{equation}\label{eq:withoutcondition}
	\Ed_{T_t,X_t,p}[\tau_\ell\wedge n]\leq \frac{\ell 2^{\ell-1}}{p^\ell}\Delta + 1.\end{equation}
	for each $0\leq t\leq n$.
	To finish, we show that $\Pd_{\mu,p}(\mathcal{G}(\Delta))\geq 1 - n^{-16}$ for an appropriate choice of $\Delta$. Indeed, by union bound and Markov's inequality:
	\[1 - \Pd_{\mu,p}(\mathcal{G}(\Delta)) \leq n\sum_{t=0}^n \Ed_{\mu,p}[\Pd_{T_t,X_t,p}(\eta_\Delta\leq n)].\]
	Now by the Simple Markov property, for $0\leq t\leq n$
	\[\Ed_{\mu,p}[\Pd_{T_t,X_t,p}(\eta_\Delta\leq n)] = \Pd_{\mu,p}(\exists \; t\leq s\leq t+n\,:\,d_{T_s}(X_s)\geq \Delta)\leq \Pd_{\mu,p}(\eta_\Delta\leq 2n).\]
	And again, by condition (\ref{eq:condicaolegal1}),
	\[1 - \Pd_{\mu,p}(\mathcal{G})\leq  n\,\Pd_{\mu,p}(\eta_\Delta\leq 2n)\leq n\,\sum_{t=0}^{2n}\Pd_{\mu,p}(d_t(X_t)\geq \Delta)\leq C\,n\,(2n+D)\,e^{-\alpha\,\Delta}.\]
	Taking \[\Delta:= \left\lceil \frac{\log(2n+D) +17\log n + \log C}{\alpha}\right\rceil\]
	guarantees $\Pd_{\mu,p}(\mathcal{G}(\Delta))\geq 1 -n ^{-16}$. Plugging this choice of $\Delta$ into (\ref{eq:withoutcondition}) gives (\ref{eq:estimateescapefinite}) inside the event $\mathcal{G}(\Delta)$.\end{proof}

\subsection{Forgetfulness of empirical measures: proof of Lemma \ref{lem:forget}}\label{ss:forget}
The proof of the lemma will follow as a consequence of another result that highlight the relation between ``escapabability" and forgetfulness. As we will see, it implies an approximated renewal property of the \textit{BGRW} process in the sense that after a random time it may be coupled to a \textit{BGRW} process starting from a semi-infinite path. Then, both walkers see around them the same tree structure as long as this coupling works. 

For this purpose, we will need additional notation. Let $Q^*$ be a semi-infinite path, i.e.,
\begin{equation}\label{def:qstar}
Q^* = \{v_0,v_1,v_2,\dots\}
\end{equation}
rooted at $v_0$ and consider a trajectory $(T^*_t,X_t^*)_{t\geq 0}$ of \textit{BGRW} started from $(Q^*,v_0)$ with empirical measures $\widehat{P}^*_n$. Also, for each  pair of integers $r$ and $n$, let $\pi(r,n)$ be
\begin{equation}
\pi(r,n) := \Pd_{Q^*_0,v_0,p}\left( \exists m \le n, \; X^*_m = x_r\right).
\end{equation}
\begin{lemma}\label{lemma:samelimit}Let $\psi:\sT_*\to\R$ be a $r$-local and bounded and $\mu$ be a $p$-escapable distribution. Consider a trajectory $(T_t,X_t)_{t\geq 0}$ of BGRW starting from $\mu$, with empirical measures~$\widehat{P}_n$. Then, 
	\begin{itemize}
		\item[(a)] for every $\varepsilon\in(0,1)$, there exist a coupling of the processes $(T_t,X_t)_{t\geq 0}$ and $(T^*_t,X^*_t)_{t\geq 0}$ such that:
		\[
		\Pd\left(\lim_{n\to +\infty}|(\widehat{P}_n - \widehat{P}^*_n)\psi|=0\right) \ge  1 - \varepsilon.
		\]
		\item[(b)] 
		\begin{equation*}
		\left | \Ed_{\mu,p} \left[ \widehat{P}_n\psi\right] - \Ed_{Q^*,v_0,p} \left[ \widehat{P}^*_n\psi\right]\right | \le 2\|\psi\|_{\infty}\pi(r,n) + \frac{2\|\psi\|_{\infty}\Ed_{\mu,p}[\tau_{2r}\wedge n]}{n+1}.
		\end{equation*}
	\end{itemize}

\end{lemma}

\begin{proof}
	\underline{Part (a).}
	Take an $\ell\in\N$ bigger than $2r$. Consider a chain started from measure $\mu$. Since~$\mu$ is $p$-escapable, the stopping time $\tau_\ell$ is finite $\Pd_{\mu,p}$-almost-surely. 
	
	Let $y_0=X_{\tau_\ell}$, $y_1$, $\dots$, $y_{\ell}$ be the unique vertices at distances $0,1,2,\dots,\ell$ at time $\tau_\ell$. Define a random time $\tau'_{\ell,r}$ as follows:
	\[\tau'_{\ell,r}:=\inf\{t\geq 0\,:\,X_{\tau_\ell + t} = y_{\ell-r}\}.\]
	Notice that this time may well be infinite. In fact, since ${\rm dist}(X_{\tau_\ell},y_{\ell-r})=\ell-r$, the second statement in Lemma \ref{lemma:returntime}, combined with the strong Markov property and the $p$-escapability of $\mu$, implies.  
	\[\Pd_{\mu,p}\left(\tau'_{\ell,r}<+\infty\right)\leq C\,e^{-\beta\,(\ell-r)}.\] 
	Note that $\tau_\ell$ is the time $X_t$ is at the tip of a path of length $\ell$ and $\tau_\ell+\tau'_{\ell,r}$ is the first time at which $X_t$ comes within distance $r$ of the other end of that path.  
	
	Now consider the process $(T^*_t,X^*_t)$ started on vertex $v_0$ of the infinite path $v_0,v_1,\dots,v_m,\dots$. Let 
	\[\tau^*_{\ell,r}:=\inf\{t\geq 0\,:\,X^*_{t} = v_{\ell-r}\}.\]
	This time may also be infinite. It is easy to see that one can couple $(T^*_t,X^*_t)$ to $(T_{\tau_\ell+t},X_{\tau_\ell+t})$ so that $\tau_{\ell,r}=\tau'_{\ell,r}$ and $[T^*_t,X^*_t]_r = [T_{\tau_\ell+t},X_{\tau_\ell+t}]_r$ for all $0\leq t\leq \tau'_{\ell,r}$. Under this coupling, if~$\tau'_{\ell,r} > n$ and $\tau_{\ell} \le n$, then 
	\begin{eqnarray}
		\widehat{P}_n\psi  &=& \frac{1}{n+1}\sum_{i=0}^n\psi(T_i,X_i) \nonumber \\ 
		\mbox{(use coupling and \textit{r}-locality)} &=&  \frac{1}{n+1}\sum_{i=0}^{\tau_\ell \wedge n-1}\psi(T_i,X_i) +  \frac{1}{n+1}\sum_{t=0}^{n-\tau_{\ell}\wedge n}\psi(T^*_{t},X^*_{t}) \nonumber \\ 
		& =& \label{eq:pstar}  \widehat{P}^*_n\psi + \frac{1}{n+1}\sum_{i=0}^{\tau_\ell \wedge n-1}(\psi(T_i,X_i) - \psi(T_{n-\tau_\ell \wedge n+i+1},X_{n-\tau_{\ell}\wedge n+i+1})).
	\end{eqnarray}
	Since $\tau_{\ell}<+\infty$ is bounded and $\widehat{P}_n\psi \le \|\psi\|_{\infty}$, this implies that, under the coupling:
	\[
	\tau'_{\ell,r}>n\Rightarrow |\widehat{P}_n\psi - \widehat{P}^*_n\psi|\leq \frac{2\|\psi\|_{\infty}\,(\tau_{\ell}\wedge n)}{n+1} \to 0.
	\]
	So, for any fixed $\ell$, we have
	\[\Pd\left(\lim_{n\to +\infty}|(\widehat{P}_n - \widehat{P}^*_n)\psi|=0\right)\geq\Pd(\tau'_{\ell,r}=+\infty)\geq  1- C\,e^{-\beta\,(\ell-r)}.\]
	In particular, this can be made greater than $ 1-\varepsilon$ with an appropriate choice of $\ell$. This proves part~$(a)$.
	
	\underline{Part (b).} Put $\ell = 2r$. From Equation (\ref{eq:pstar}) and the fact that under the coupling $\tau'_{\ell,r} = \tau^*_{\ell,r}$, we may obtain that
	\begin{equation}
	\begin{split}
	\left | \widehat{P}_n\psi - \widehat{P}^*_n\psi\right | & = \left | \widehat{P}_n\psi - \widehat{P}^*_n\psi\right |\mathbb{1}\{\tau'_{\ell,r} \le n \} + \left | \widehat{P}_n\psi - \widehat{P}^*_n\psi\right |\mathbb{1}\{\tau'_{\ell,r} > n \} \\
	& \le 2 \|\psi \|_{\infty}\mathbb{1}\{\tau^*_{\ell,r} \le n \} + \frac{2\|\psi\|_{\infty}\,(\tau_{\ell}\wedge n)}{n+1}
	\end{split}
	\end{equation}
	taking the expected value both sides it is enough to proves part (b) and finally proves the lemma.
\end{proof}

Now we see how Lemma \ref{lem:forget} follows immediately from the above result.
\begin{proof}[Proof of Lemma \ref{lem:forget}] 
%
	
	From part (a) of Lemma~\ref{lemma:samelimit}, one sees that, given a deterministic $c\in\R$,
	\[\Pd_{\mu,p}(\widehat{P}_n\psi\to c) = 1\Leftrightarrow \Pd_{Q^*,v_0,p}(\widehat{P}^*_n\psi\to c) = 1.\]
	This holds for {\em any} $p$-escapable measure $\mu$. In particular, if $\nu$ is {\em also} $p$-escapable $\widehat{P}_n\psi\to c$ $\Pd_{\nu,p}$-a.s., the same must hold with $\nu$ replacing $\mu$. In other words, Lemma~\ref{lemma:samelimit} (a) implies the lemma.

\end{proof}

\subsection{Convergence of average of local functions: proof of Lemma \ref{lema:rfunclema}}\label{ss:rfunclema}
The proof of Lemma~\ref{lema:rfunclema} will follow essentially from the lemma below that states that, when started from a finite initial condition, the average of local functions is close to its expected value but started from the semi-finite path $Q^*$ introduced in (\ref{def:qstar}).
\begin{lemma}\label{lemma:ME}
	Given $p_0>0$, there exist constants $C,\alpha>0$ depending only on $p_0$ such that, for all  $p_0\leq p\leq 1$, all finite $(T_0,x_0)$, all $-1\le \psi \le 1$  $r$-local, all $\varepsilon >0$,  and all $1\leq q<n$, $\lfloor n/q \rfloor>0$,
	\begin{equation}
	\Pd_{T_0,x_0,p} \left( \left | \widehat{P}_n\psi - \Ed_{Q^*,v_0,p}\left[ \widehat{P}_q\psi \right] \right | \ge \varepsilon + \frac{2q}{n} + err(q,r) \right) \le 2\exp\{ - \varepsilon^2\lfloor n/q \rfloor /32\} + \frac{1}{q^{16}},
	\end{equation}
	where 
	\begin{equation}
	err(q,r) :=  4\pi(r,q) +
	\frac{r 2^{2r}}{q\, p^{2r}}\left\lceil \frac{\log(2q+|V(T_0)|) + 17\log q +\alpha \Delta(T_0)\log C }{\alpha} \right\rceil + \frac{1}{q}.
	\end{equation}
\end{lemma}

\begin{proof}
	Let $n=kq+s$, with $k\geq 1$. Obseve that, for every $k,q\geq 1$, we can write 
	\begin{equation}
	\begin{split}
	\widehat{P}_{kq-1}\psi= \frac{1}{kq}\sum_{i=0}^{kq-1}\psi(T_i,X_i) &= \frac{1}{k}\left(\frac{1}{q}\sum_{i=0}^{q-1}\psi(T_i,X_i) + \cdots + \frac{1}{q}\sum_{i=kq-q}^{kq-1}\psi(T_i,X_i)\right)
	\\
	&= \frac{1}{k} \sum_{j=1}^k S_j\;,
	\end{split}
	\end{equation}
	where,   $S_j := \displaystyle \frac{1}{q}\sum_{i=(j-1)q}^{jq-1}\psi(T_i,X_i)$, for $j \in \{1, \cdots, k\}$. 
	Let us denote by $\{\tilde{\mathcal{F}}_j\}_{j \ge 0}$ the following filtration:
	\begin{equation}
	\tilde{\mathcal{F}}_0 = \mathcal{F}_0; \; \; \tilde{\mathcal{F}}_j = \mathcal{F}_{jq-1}.
	\end{equation}
	By the definition of $S_j$, the process $\{S_j\}_{j \ge 1}$ is adapted to the filtration~$\{\tilde{\mathcal{F}}_j\}_{j \ge 0}$. Moreover, the process $\{M_k\}_{k \ge 1}$ defined as
	\begin{equation}
	M_k := \sum_{j=1}^{k} S_j - \Ed_{T_0,x_0,p}\left[ S_j \middle | \tilde{\mathcal{F}}_{j-1}\right] = M_{k-1} + S_k - \Ed_{T_0,x_0,p}\left[ S_k \middle | \tilde{\mathcal{F}}_{k-1}\right]
	\end{equation}
	is a mean zero martingale with respect to $\{\tilde{\mathcal{F}}_j\}_{j \ge 0}$, whose increments are bounded by $2$, since~$S_j$ is an average of random variables bounded by $1$.

	With all these definitions, by triangle inequality, we have that 
	\begin{equation}\label{ineq:triangle}
	\begin{split}
	\left | \widehat{P}_n\psi - \Ed_{Q^*,v_0,p}\left[ \widehat{P}_q\psi \right] \right | 
	& \le
	\left | \widehat{P}_n\psi -  \widehat{P}_{kq-1}\psi  \right | \quad (a)
	\\		
	&		+
	\left | \widehat{P}_{kq-1}(\psi) - \frac{1}{k}\sum_{j=1}^k\Ed_{T_0,x_0,p}\left[ S_j \middle | \tilde{\mathcal{F}}_{j-1}\right] \right  | \quad (b)
	\\
	&
	+
	\left | \frac{1}{k}\left(\sum_{j=1}^k \Ed_{T_0,x_0,p}\left[ S_j \middle | \tilde{\mathcal{F}}_{j-1}\right] - \Ed_{Q^*,v_0,p}\left[ \widehat{P}_q\psi \right] \right)  \right | \quad (c).
	%
	\end{split}
	\end{equation}
	
	As far as $(a)$  is concerned, we observe that the following deterministic bound holds:
	\begin{equation*}
	\begin{split}
	\left | \widehat{P}_n\psi -
	\widehat{P}_{kq-1}\psi \right |
	& = \left | \left(\frac{kq}{n+1}-1\right)\frac{1}{kq}\sum_{i=0}^{kq-1}\psi(T_i,X_i) + \frac{1}{n+1}\sum_{i=kq	}^{kq+s}\psi(T_i,X_i)\right | \\
	& \le \frac{2 \|\psi \|_{\infty}(s+1)}{n+1} < \frac{2q}{n},
	\end{split}
	\end{equation*}
	since, by hypothesis, $\|\psi \|_{\infty} \le 1$. Thus, to prove the claim, it suffices to show that
	\begin{equation}\label{eq:pippo}
	\Pd_{T_0,x_0,p} \Big( (b) + (c) \ge \varepsilon + err(q,r) \Big) 
	\le 2\exp\{ - \varepsilon^2k/32\} + \frac{1}{q^{16}}\;.	
	\end{equation}
	Using the simple observation that
	\[
	\Pd_{T_0,x_0,p} \Big( (b) + (c) \ge \varepsilon + err(q,r) \Big)\leq 
	\Pd_{T_0,x_0,p} \left( (b) \ge \frac{\varepsilon + err(q,r)}{2} \right) + 	\Pd_{T_0,x_0,p} \left( (c) \ge \frac{\varepsilon + err(q,r)}{2} \right)\;,
	\]
	we will prove Equation~\eqref{eq:pippo} by showing 
	\begin{align*}
	&(1): \qquad 	\Pd_{T_0,x_0,p} \left( (b) \ge \frac{\varepsilon }{2} \right)\le 2\exp\{ - \varepsilon^2k/32\}\,
	\\
	&(2): 
	\qquad 	\Pd_{T_0,x_0,p} \left( (c) 
	> 
	\frac{ err(q,r)}{2} \right) \le 
	\frac{1}{q^{16}}\, 
	\end{align*}	

	To show $(1)$, we observe that the martingale  $M_k$ can be written  as the following difference
	\begin{equation}
	M_k = k\widehat{P}_{kq-1}(\psi)  -\sum_{j=1}^k\Ed_{T_0,x_0,p}\left[ S_j \middle | \tilde{\mathcal{F}}_{j-1}\right].
	\end{equation}
	Thus, $(b)$ is equal to $|M_k/k|$ and we apply Azuma's inequality which gives us:
	\begin{equation}\label{eq:Mk}
	\Pd_{T_0,x_0,p}\left( \left|M_k\right| \ge \frac{\varepsilon k}{2}\right) \le 2\exp\{ - \varepsilon^2k/32\}\, .
	\end{equation}

	To show $(2)$, we first observe that, by the Simple Markov Property, for all $j \in \{2, \cdots, k\}$ the following identity holds  
	\begin{equation}
	\Ed_{T_0,x_0,p}\left[ S_j \middle | \tilde{\mathcal{F}}_{j-1}\right] = \frac{1}{q}\sum_{i=1}^{q}\Ed_{T_{(j-1)q-1},X_{(j-1)q-1},p} \left[ \psi(T_i,X_i)\right] = \Ed_{T_{(j-1)q-1},X_{(j-1)q-1},p} \left[ \widehat{P}_q(\psi)\right].
	\end{equation}
	Then, by Part (b) of Lemma \ref{lemma:samelimit} we have
	\begin{equation}\label{ineq:SMQ}
	\begin{split}
	\left| \Ed_{T_0,x_0,p}\left[ S_j \middle | \tilde{\mathcal{F}}_{j-1}\right] - \Ed_{Q^*,v_0,p}\left[ \widehat{P}_q(\psi) \right]\right| & \le 2\pi(r,q) + \frac{2\Ed_{T_{(j-1)q-1},X_{(j-1)q-1},p} \left[\tau_{2r}\wedge q\right]}{q}\, ,
	\end{split}
	\end{equation}
	and implies
	\[
	(c) \le 2\pi(r,q) + \frac{1}{k}\sum_{j=1}^k\frac{2\Ed_{T_{(j-1)q-1},X_{(j-1)q-1},p} \left[\tau_{2r}\wedge q\right]}{q}\;.
	\]
	Therefore,  it is enough to show that the following events 
	\begin{equation}\label{eq:baggio}
	B:=\left\{ \sum_{j=0}^{k-1}\frac{2\Ed_{T_{jq-1},X_{jq-1},p} \left[\tau_{2r}\wedge q\right]}{k\,q} > \frac{r 2^{2r}}{q\, p^{2r}}\left\lceil \frac{\log(2q+|V(T_0)|) + 17\log q +\alpha \Delta(T_0)\log C }{\alpha} \right\rceil + \frac{1}{q} \right\}
	\end{equation}
	happens with probability at most $q^{-16}$. 
	
	By Lemma \ref{lem:quantitative} (using   the constants coming from Lemma~\ref{lem:degreetail}) we know that, if we define
	\[
	A_t:=\left\{ \Ed_{T_t,X_t,p}[\tau_{2r}\wedge q] \leq \frac{r 2^{2r}}{p^{2r}}\left\lceil \frac{\log(2q+|V(T_0)|) + 17\log q }{\alpha} + \Delta(T_0)\log C\right\rceil+ 1\right\}\, ,
	\]
	it holds that 
	\[
	\Pd_{T_0,x_0,p}\left( \bigcap_{t=0}^q A_t\right) \geq  1-\frac{1}{q^{16}}\;, 
	\]
	which implies $
	\Pd_{T_0,x_0,p}\left( B \right) \leq \frac{1}{q^{16}}$, and concludes the proof.
\end{proof}

\begin{corollary}\label{cor:totti} Given  $p_0>0$, for all $p\in [p_0, 1]$, for all finite $(T_0,x_0)$ and for all $r$-local function $\psi$  such that $-1\leq \psi\leq 1$,    $\{P_n\psi\}_{n\in \N}$ is  a Cauchy sequence $\mathbb{P}_{T_0,x_0}$-a.s.
\end{corollary}

\begin{proof} The proof will follow from  the claim below.
	
	\begin{claim} Given $p_0>0$ and $(T_0, x_0)$ finite, there exists $n_0$ and a constant $c$ depending only on $p_0, T_0$ such that, for all $n\geq n_0$ and for all $p\in [p_0,1]$: 
		\begin{equation}
		\Pd_{T_0,x_0,p}\left( \left| \widehat{P}_m\psi - \widehat{P}_{m'}\psi\right|> \frac{c}{\log n}, \text{ for some }m,m'\in \{n,n+1,\cdots,n^2\}\right) \le \frac{4}{n^{4}}.
		\end{equation}
	\end{claim}	
	
	\begin{claimproof}
		To prove the claim we show that there	exists $n_0$ such that for all $n\geq n_0$ and for all  $m, m'\in \{n,n+1,\cdots,n^2\}$, 
		\begin{equation}\label{eq:minnie0}
		\Pd_{T_0,x_0,p}\left( \left| \widehat{P}_m\psi -\widehat{P}_{m'}\psi  \right | > \frac{c}{\log n}\right)  \le \frac{4}{n^{8}}\ ,
		\end{equation}
		and then use union bound. We proceed by showing that for all  $m\in \{n,n+1,\cdots,n^2\}$ we can find values for $r$ and $q$ such that
		\begin{equation}\label{eq:minnie}
		\Pd_{T_0,x_0,p}\left( \left| \widehat{P}_m\psi - \Ed_{Q^*,v_0}\left[ \widehat{P}_q(\psi) \right] \right | >\frac{c}{2\log n} \right)  \le \frac{2}{n^{8}}\ ,
		\end{equation}
		which implies Equation~\eqref{eq:minnie0} by using triangle inequality. In order to show that Equation~\eqref{eq:minnie} holds, 	we use Lemma~\ref{lemma:ME}, suitably choosing the values
		of $q$, $r$ and $\varepsilon$. 
		We begin observing that by Corollary~\ref{cor:etaybound} we have that 
		\[
		\pi(r,q)\leq \frac{C_1}{\log n}\, .
		\]			
		Choosing $q= n^{1/2}$ and $r= \left \lfloor \frac{\log n}{8\log(2/p)} \right \rfloor$, we obtain that 
		\begin{equation}\label{eq:bounderr}
		\begin{split}
		err(q,r) & =  4\pi(r,q) +
		\frac{4r 2^{2r}}{q\,p^{2r}}\left\lceil \frac{\log(2q+|V(T_0)|) + 17\log q +\alpha \Delta(T_0)\log C }{\alpha} \right\rceil + \frac{1}{q} \\
		& \le \frac{4C_1}{\log n}+ \frac{\log n}{2 n^{1/4}\log(2/p)} \left( \frac{ \log 2n^{1/2} + |V(T_0)| + 9\log n}{\alpha} + \Delta(T_0)\log C \right) +\frac{1}{n^{1/2}} \\
		& 
		\le \frac{4C_1}{\log n}+ \frac{\log^2 2n}{\alpha \,n^{1/4} } + \frac{|V(T_0)|\log n }{\alpha n^{1/4}} + \frac{9\log^2 n }{\alpha n^{1/4}} +  \frac{\log n \,\Delta(T_0) \log C}{n^{1/4}} + \frac{1}{n^{1/2}}
		\end{split}
		\end{equation}
		Given that the constants $C, C_1, \alpha$ only depend on $p_0$, and $T_0$ is finite, there exists a constant $c$ and a sufficiently large $n_0$, both depending on $p_0$ and $T_0$ only,   such that 
		\[
		err(q,r)\leq \frac{c}{6\log n}\;.
		\]		
		Setting $\e= c/6\log^{-1} n$ and since $2q/m\leq  2n^{-1/2}\leq c/6\log^{-1} n$ for sufficiently large $n$, 
		by Lemma~\ref{lemma:ME},  for all $m \in \{n,n+1,\cdots,n^2\}$, we have that 
		\begin{equation}\label{eq:me2}
		\Pd_{T_0,x_0,p} \left( \left | \widehat{P}_m\psi - \Ed_{Q^*,v_0}\left[ \widehat{P}_q\psi\right] \right | > \frac{c}{2\log n}\right) \le 2e^{ - \varepsilon^2/32 \sqrt{n}} + \frac{1}{q^{16}} \le \frac{2}{n^8}\, ,
		\end{equation}
		for sufficiently large $n$ (how large depends on $\varepsilon$ which only depends on $p_0, T_0$). 		
		Thus, for all ~$m, m'\in \{n,n+1,\cdots,n^2\}$ we have 
		\begin{equation}
		\Pd_{T_0,x_0,p} \left( \left | \widehat{P}_m\psi - \widehat{P}_{m'}\psi \right | >  \frac{c}{\log n}\right)  \le \frac{4}{n^{8}}\ ,
		\end{equation}
		and finally, by union bound,  the claim follows.

	\end{claimproof}

	The above claim, together with a Borel-Cantelli argument, implies Corollary~\ref{cor:totti}, i.e,  $\mathbb{P}_{T_0,x_0,p}$ almost surely there exist~$n_0$ such that, for all $n\ge n_0$ and  for all $m,m' \in \{n,n+1, \cdots, n^2\}$ we have
	\[
	\left | \widehat{P}_m\psi - \widehat{P}_{m'}\psi \right | \le \frac{c}{\log n}.
	\]
	For any $\delta >0$, we can choose $n_1$ large enough so that $\frac{c}{\log n} < \delta/2$  and $n_1 \ge n_0$. Now, for~$m, m' \ge n_1$, assuming $m \le m'$, there exists $k_0$ such that $m' \in \{m^{2^{k_0}},\cdots, m^{2^{k_0+1}}\}$. Thus, 
	\begin{equation}
	\begin{split}
	\left | \widehat{P}_m\psi - \widehat{P}_{m'}\psi \right | & \le \sum_{j=0}^{k_0-1} \left | \widehat{P}_{m^{2^j}}\psi - \widehat{P}_{m^{2^{j+1}}}\psi \right | + \left | \widehat{P}_{m^{2^{k_0}}}\psi - \widehat{P}_{m'}\psi \right | \\ 
	& \le \sum_{j=0}^{k_0}\frac{c}{({2^j \log m})} \le  \frac{2c}{\log m} \le \frac{2c}{\log n_1}  < \delta\, .
	\end{split}
	\end{equation}
\end{proof}
Now we are able to prove Lemma~\ref{lema:rfunclema}.
\begin{proof}[Proof of Lemma~\ref{lema:rfunclema}]We prove first that the sequence $\{\widehat{P}_n\psi\}_{n \in \N }$ started from any finite initial condition always converges to the same limit. Finally, we prove that this limit must be a constant. 
	
	Consider two independent \textit{BGRW} processes $(T_t,X_t)_{t \in \N}$ and $(T'_t,X'_t)_{t \in \N}$ starting from two finite initial conditions $(T_0,x_0)$ and $(T'_0,x_0')$ respectively. Let $\{\widehat{P}_n\psi\}_{n \in \N }$ and $\{\widehat{P}'_n\psi\}_{n \in \N }$ be the averages associated to each process. The proof follows from (\ref{eq:me2}) and triangle inequality, which combined give us that
	\begin{equation}
		\begin{split}
		\Pd \left( 	\left|\widehat{P}_n\psi - \widehat{P}_{n}'\psi\right| \ge \frac{c}{\log n} \right)  \le \frac{4}{n^8}.
		\end{split}
	\end{equation}
	for some positive constant $c$ depending on $p$, $T_0$ and $T_0'$ only. Then,	
 an application of Borel-Cantelli Lemma allows us to conclude that
	both processes converge to the same limit $\mathcal{M}_p\psi$. Once we have that, we finally observe that since the \textit{BGRW} processes are independent, the events $\{\mathcal{M}_p\psi \in A\}$ and~$\{\mathcal{M}_p\psi \in B\}$ are independents for all Borel sets on the real line. This implies that $\mathcal{M}_p\psi$ is independent of itself, thus, constant almost surely.

\end{proof}





\section{Weak convergence of empirical measures}\label{sec:localweak}

In this section we prove  Theorem~\ref{thm:localweak} and Theorem~\ref{thm:oneended} in the Introduction. In fact, we prove slightly stronger statements that only require  the initial measure to be $p$-escapable, which in particular applies to all initial measures supported on finite trees (cf. Lemma \ref{lem:quantitative}).
Specifically, we will show that, for any $p$-escapable initial measure $\mu$, there exits a measure~$P_p$ on the space $\sT_*$ such that,  when $n\to +\infty$, 
\begin{align*}
\widehat{P}_n&:=\frac{1}{n+1}\sum_{t=0}^{n}\delta_{[T_t,X_t]}\Rightarrow P_p\, \quad \mbox{ $\Pd_{\mu,p}$-a.s.}
\\
P_n &:=\frac{1}{n+1}\sum_{t=0}^n\,\mu K_p^t\Rightarrow P_p\,
\end{align*}
where  $K_p$  denotes the Markov transition kernel of the BGRW process over $\sT_*$.

In our proof we use  the next proposition, which gives a convenient criterion for checking weak convergence of distributions over $\sT_*$. Before stating the proposition we need to define a specific kind of local function.
Given $r\in\N$ and~$[T,o]\in\sT_*$, let $D_r([T,o])$ denote the largest degree of a vertex in $T$ at distance $r$ from $o$:
\[D_r(T,o):=\max\{d_T(x)\,:\, x\in V(T),\, {\rm dist}_T(x,o)\leq r\}.\]
One can check that $D_r:\sT_*\to\R$ is $(r+1)$-local. 

\begin{proposition}[Proof in Appendix \ref{sec:criterionlocalweak}]\label{prop:convergencelocal} Suppose $\{Q_n\}_{n\in\N}$ is a sequence of probability measures over $\sT_*$. Then:
	\begin{enumerate}
		\item {\em Tightness:} If $\lim_{k\to +\infty}\limsup_{n\in\N}Q_n(\{[T,o]\in\sT_*\,:\, D_r(T,o)\geq k\})=0$ for any $r\geq 0$, then $\{Q_n\}_{n\in\N}$ is tight.
		\item {\em Weak convergence:} Assume $\{Q_n\}_{n\in\N}$ satisfies item (1). Then there exists a countable family $\mathcal{S}$ of bounded local functions such that, if $\lim_n Q_nf$ exists for all $f\in\mathcal{S}$, then $Q_n\Rightarrow Q$ for some $Q$ that is uniquely defined by:
		\[Qf = \lim_n Q_nf,\, f\in\mathcal{S}.\]  
	\end{enumerate}
\end{proposition}
%
%
%
We begin with addressing tightness for $\{P_n\}_{n\in \N}$. The next Lemma shows C\`{e}saro-mean-type sequences of measures generated by our Markov chain are tight when the initial measure satisfy the quantitative criterion of $p$-escapability defined in~Lemma~\ref{lem:quantitative}. 
\begin{lemma}\label{lem:tightness} Assume $\mu$ is a measure over $\sT_*$ with the following property: there exist $C,\alpha,D>0$ such that for all $n,k\geq 1$:
	\begin{equation}\label{eq:condicaolegal}
	\sum_{t=0}^n\Pd_{\mu,p}(\degr{t}{X_t}\geq k)\leq C\,(n+D)\,e^{-\alpha\,k}.
	\end{equation}
	Then the sequence
	\[P_n:= \frac{1}{n+1}\sum_{t=0}^n\mu K_p^t\]
	satisfies the tightness criterion in Proposition \ref{prop:convergencelocal}. In fact, the following quantitative estimate is satisfied: for all $r,k\geq 1,n\in\N$:
	\[P_n(\{[G,o]\in\sT_*\,:\, D_r(G,o)\geq k\})\leq C\,\frac{(n+D+r)}{n+1}\,\left(\frac{(k+1)^{r+1}-1}{k}\right)\,e^{-\alpha\,k}.\]\end{lemma}
\begin{proof}Note that:
	\[(n+1)P_n(\{[G,o]\in\sT_*\,:\, D_r(G,o)\geq k\}) = \sum_{t=0}^n\Pd_{\mu,p}(D_r([T_t,X_t])\geq k).\]
	So it suffices to prove the following quantitative estimate: for all $r,n,k\geq 1$,
	\begin{equation}\mbox{\bf Goal: } \sum_{t=0}^n\Pd_{\mu,p}(D_r(T_t,X_t)\geq k)\leq C\,(n+D+r)\,\left(\frac{(k+1)^{r+1}-1}{k}\right)\,e^{-\alpha\,k},\end{equation}
	which is true for $r=0$ by our assumption (\ref{eq:condicaolegal}).
	
	Consider $r>0$ and assume that we have proven that, for each $j\leq r-1$, there exists $C_j>0$ depending only on $j$ and $p$, and an exponent $\alpha$ as in the base case, such that
	\[\mbox{\bf (induction hyp.) } \forall m\in\N\,:\,\sum_{t=0}^m\,\Pd_{\mu,p}(D_{j}(T_t,X_t)\geq k)\leq C_{j}\,(m+D+j)\,e^{-\alpha\,k},\]
	for all $k\geq 1$ (again, this is true for $j=0$ with $C_0=C$). We {\em claim} that a similar statement holds for $j=r$ with $C_r = C_0 + kC_{r-1}$. A simple induction would then imply our goal:
	\[C_r=C_0\,\sum_{i=0}^r(k+1)^i = \frac{(k+1)^{r+1}-1}{k}\,C.\]
	To obtain our bound, note that 
	\begin{eqnarray*}\sum_{t=0}^n\,\Pd_{\mu,p}(D_{r}(T_t,X_t)\geq k)&\leq & \sum_{t=0}^n\,\Pd_{\mu,p}(\degr{t}{X_t}\geq k) \\ & & + \sum_{t=0}^n\,\Pd_{\mu,p}(D_{r}(T_t,X_t)\geq k,\degr{t}{X_t}\leq k).\end{eqnarray*}
	The first term in the RHS is $\leq C_0\,(n+D)\,e^{-\alpha\,k}$ by the base case. For the second term, we observe that $D_{r}(T_t,X_t)\geq k$ -- i.e. some vertex at distance $r$ from $X_t$ has degree $\geq k$ -- means $D_{r-1}(T_t,v)\geq k$ for one of the neighbors of $X_t$ in $T_t$. Now, if $\degr{t}{X_t}\leq k$, there is a chance of at least $ 1/(k+1)$ that $X_{t+1}=v$ and thus $D_{r-1}(T_{t+1},X_{t+1})\geq k$. We conclude:
	\[\frac{\Pd_{\mu,p}(D_{r}(T_t,X_t)\geq k,\degr{t}{X_t}\leq k)}{k+1}\leq \Pd_{\mu,p}(D_{r-1}(T_{t+1},X_{t+1})\geq k).\]
	So: \begin{eqnarray*}\sum_{t=0}^n\,\Pd_{\mu,p}(D_{r}(T_t,X_t)\geq k) &\leq & C_0\,(n+D)\,e^{-\alpha\,k}\\ & & + (k+1)\,\sum_{t=0}^n\,\Pd_{\mu,p}(D_{r-1}(T_{t+1},X_{t+1})\geq k) \\ \mbox{(induction hyp. for $j=r-1,m=n+1$)} 
		&\leq & C_0\,(n+D)\,e^{-\alpha\,k} \\ & & + (k+1)C_{r-1}\,(n+r+D) \,e^{-\alpha k}\\&\leq &C_r\,(n+r+D)\,e^{-\alpha\,k},\end{eqnarray*}
	where $C_r:= C + (k+1)\,C_{r-1}$ as desired.\end{proof}
In the next result -- a strengthening of Theorem \ref{thm:localweak} --, we prove convergence of empirical measures $\widehat{P}_n$ and Ces\`{a}ro-style convergence of the chain $K_p$ from any initial measure that is $p$-escapable. We also characterize the limiting probability measure of the chain as a stationary distribution for $K_p$.

\begin{theorem}[Convergence of the empirical measure; proof in \S \ref{sec:proof.convergence}]\label{thm:convergence} Given $p\in (0,1]$, there exists a probability measure $P_p$ on the space $\sT_*$ of rooted trees such that, for any $p$-escapable initial measure $\mu$, \[\widehat{P}_n\Rightarrow P_p\;\mbox{ $\Pd_{\mu,p}$-a.s. \quad   and \quad }\frac{1}{n+1}\sum_{t=0}^n\,\mu K_p^t\Rightarrow P_p \;\mbox{ when $n\to +\infty$}.\] For $0<p_0\leq p\leq 1$, the measure $P_p$ satisfies:
\begin{equation}\label{eq:escapablestationary}\forall k,r\geq 1\,:\, P_p(\{[T,o]\in\sT\,:\,D_r([T,o])\geq k\})\leq C\,\frac{(k+1)^{r+1}-1}{k}\,e^{\alpha\,(2-k)}\end{equation}
where $C,\alpha>0$ only depend on $p_0$. Moreover, $P_p$ is an invariant measure for $K_p$.\end{theorem}

Note that the measure $P_p$ itself is $p$-escapable, by Lemma \ref{lem:quantitative} combined with (\ref{eq:escapablestationary}). That is,~$P_p$ is the unique $p$-escapable invariant measure for $K_p$. On the other hand, the hard tree in Example \ref{ex:hardtree} shows that one does not have convergence to $P_p$ from all initial measures. 

Our next result extends Theorem \ref{thm:oneended} from the introduction.

\begin{theorem}[Support of $P_p$; \S \ref{sec:proof.oneended2}]\label{thm:oneended2} Given $0<p\leq 1$,  let $P_p$ be the limit measure in Theorem \ref{thm:convergence}. Then $P_p$ is supported on one-ended infinite trees. Moreover, given a finite rooted tree $(S_0,x_0)$ of diameter $\leq h$,
\[P_p(\{[T,o]\in\sT_*\,:\, [T,o]_h = [S_0,x_0])\})>0.\]\end{theorem}

In particular, this shows that, if $\mu$ is such that $\widehat{P}_n\Rightarrow P_p$ $\Pd_{\mu,p}$-a.s., then $\Pd_{\mu,p}(\tau_\ell<+\infty) = 1$ i.e. $\mu$ is $p$-escapable. In this sense, $p$-escapability is a necessary condition in Theorem \ref{thm:convergence}.

\subsection{Weak convergence of empirical measures: proof of Theorem~\ref{thm:convergence}}\label{sec:proof.convergence}

We present here the proof of convergence to the stationary measure. 

\begin{proof}[Proof of Theorem \ref{thm:convergence}] Fix $\mu$ and $p$. Theorem \ref{thm:rfuncthm} guarantees that 
\begin{equation}\label{eq:boundedlocalconvergence}\Pd_{\mu,p}(\widehat{P}_n\psi\to \mathcal{M}_p\psi)=1\end{equation}
for all bounded local functions, where $\mathcal{M}_p\psi$ does not depend on $\mu$. If we write:
\[P_n := \frac{1}{n+1}\sum_{t=0}^n\,\mu K_p^t,\]
we also have 
\begin{equation}\label{eq:boundedlocalconvergence2}P_n\psi = \Ed_{\mu,p}[\widehat{P}_n\psi]\to \mathcal{M}_p\psi\end{equation}
by the Bounded Convergence Theorem. 

We claim that the sequences $\{P_n\}_{n}$ and $\{\widehat{P}_n\}_{n}$ are  {\em tight} (almost surely, in the second case). More specifically, we will apply the criterion in Proposition  \ref{prop:convergencelocal}, part $(1)$. In fact, fixing $k,r\in\N$, we may take
\[\phi_{k,r}(T,x):=\Ind{\{D_r(T,x)\geq k\}}\,\,([T,x]\in\sT_*)\]
which is a $(r+1)$-local function, so that
\[P_n\phi_{k,r}\to \mathcal{M}_p\phi_{k,r},\; \text{ and }\;\widehat{P}_n\phi_{k,r}\to \mathcal{M}_p\phi_{k,r}\, \;\; \mbox{ $\Pd_{\mu,p}$-a.s.}\]
Now, to evaluate the limit in the RHS, we replace the initial measure $\mu$ with a measure~$\mu_0$ supported on a finite tree with two vertices. Applying Lemma \ref{lem:tightness} in conjunction with Lemma~\ref{lem:degreetail}, and letting $n\to +\infty$, we obtain:
\[
\mathcal{M}_p\phi_{k,r}\leq C\left(\frac{(k+1)^{r+1}-1}{k}\right)\,e^{\alpha\,(2-k)}.
\]
Since the RHS of this inequality goes to $0$ as $k\to +\infty$, we obtain the required tightness. 

We now prove that the weak convergence criterion of Proposition \ref{prop:convergencelocal} (part $(2)$) is also satisfied, which asks for the convergence of expectations of all bounded local functions. For $P_n$ this is immediate from (\ref{eq:boundedlocalconvergence2}). For $\widehat{P}_n$ there is the slight issue that we have:
\[\forall \mbox{ bounded local }\psi\,:\, \Pd_{\mu,p}(\widehat{P}_n\psi\to\mathcal{M}_p(\psi))=1\]
and we now need to ``move the quantifier $\forall$ inside the probability". However, Proposition \ref{prop:convergencelocal} shows that we only need to worry about a countable family of $\psi$, and there is no problem in moving the quantifier inside for a countable family.

The upshot is that, assuming this claim, we have that Proposition \ref{prop:convergencelocal} assures that $P_n$ and $\widehat{P}_n$ converge weakly (almost surely, in the second case) to the {\em same} probability measure $P_p$, which is uniquely characterized by the fact that $P_p\psi = \mathcal{M}_p\psi$. In particular, this measure satisfies the property that:
\[\forall r,k\geq 1\,:\,P_p(\{[T,o]\,:\,D_r(T,o)\geq k)\leq C\left(\frac{(k+1)^{r+1}-1}{k}\right)\,e^{\alpha\,(2-k)}.\]

Finally, we show that $P_p$ is an invariant measure for $K_p$. 
To this aim, we need the following fact  assuring that  if $\psi$ is bounded Lipschitz, so is $K_p(\psi)$.
\begin{fact}\label{fact:bl}Given a bounded measurable function $\psi:\sT_*\to \R$, let 
\[K_p(\psi)([T,o]):=\int_{\sT_*}\,\psi([T',o'])\,K_p([T,o],{\rm d}[T',o']).\]
If $\psi$ is bounded and Lipschitz, so is $K_p(\psi)$.\end{fact}
\begin{proof}[Proof sketch] This follows from the fact that, if $[T,o]$ and $[S,v]$ satisfy $\rho([T,o],[S,v])\leq 1/(1+r)$ for some $r>0$, one can couple $[T',o']\sim K([T,o],\cdot)$ and $[S',v']\sim K([S,v],\cdot)$ so that $\rho([T',o'],[S',v'])\leq 1/r$.\end{proof}

So we may apply $P_p$ to the function $K_p(\psi)$ and obtain:
\[P_pK_p(\psi) = \lim_{n\to +\infty} P_nK_p(\psi) = \lim_{n\to +\infty} \frac{1}{n+1}\sum_{t=0}^n\mu K_ p^{t+1}(\psi) = \lim_{n\to +\infty}\frac{1}{n+1}\sum_{t=1}^{n+1}\mu K_ p^{t}(\psi) = P_p(\psi).\]
Since $\psi$ is an arbitrary bounded Lipschitz function, this implies $P_pK_p=P_p$. So $P_p$ is a stationary measure for the chain $K_p$.\end{proof}

\subsection{On the support of the limiting distribution: proof of Theorem \ref{thm:oneended2}}\label{sec:proof.oneended2}

\begin{proof}[Proof of Theorem \ref{thm:oneended2}] That $P_p$ is supported on infinite trees is immediate. We prove next the final statement in the theorem, namely: 
\[P_p(\{[T,o]\in\sT_*\,:\, [T,o]_h = [S,x_0])\})>0.\]

To see this, recall that $(S,x_0)$ has height $h$. Let $x_1$ be a leaf node of $S$ at distance $h$ from~$x_0$. One may traverse the tree $S$ by first walking $h$ steps from $x_1$ to $x_0$ and then doing a depth-first search on $S$ from $x_0$, which will end back at $x_0$. 

This traversal gives us a recipe for ``creating" $S$ inside $T_{t_0}$, for a certain $t_0$. Namely, all we need is that:
\begin{enumerate}\item The walk creates leaves and jumps to them for $h$ consecutive steps (this corresponds to moving from $x_1$ to $x_0$
\item The walk then follows the sequence of steps given by the DFS in $S$, adding leaves as needed.
\end{enumerate}
The above sequence of steps has positive probability and, if followed, ensures $[T_{t_0},X_{t_0}]_h = [S,x_0]$. By stationarity of the process
\[P_p(\{[T,o]\in\sT_*\,:\, [T,o]_h = [S,x_0])\}) = \Pd_{P_p,p}([T_{t_0},X_{t_0}]_h = [S,x_0])>0.\]

We now prove one-endedness. Denote by $O(s)$ the set of all rooted trees $[T,o]$ that satisfy the following two properties:
\begin{enumerate}
\item $T$ is infinite (but locally finite, as all other trees in $\sT_*$);
\item there exists a $h\geq s$ such that all paths of length $\geq h$ starting from $o$ pass through a single vertex at distance $s$  from $o$.
\end{enumerate}
One may check that the set of one-ended trees is precisely $\cap_{s\in\N}O(s)$. Therefore, it is enough to  show that $P_p(O(s))=1$ for all $s\in\N$.

So fix $s\in\N$ and take a trajectory $[T_t,X_t]$ of BGRW started from $P_p$. Since $P_p$ is stationary,
\[\forall n\in\N\,:\, \Pd_{P_p,p}([T_n,X_n]\in O(s)) = P_p(O(s)).\]
We will show that, for all $\varepsilon>0$, one can choose $n\in\N$ so that the LHS above is $\geq 1-\varepsilon$, which suffices to finish the proof. 

Recall that $P_p$ is $p$-escapable, as we commented after the statement of Theorem \ref{thm:convergence}. This means that $\tau_{\ell}<+\infty$ almost surely for all $\ell\in\N$. Take $\ell\geq s$. We will follow the reasoning in the proof Lemma \ref{lem:forget} -- more specifically the Claim at the end of the proof -- and argue as follows. Let $X_{\tau_\ell}=v_0,v_1,\dots,v_{\ell}$ is the (induced) path of length $\ell$ around $X_{\tau_\ell}$ at time $\tau_\ell$. Define
\[\tau'_{\ell,s}:=\inf\{t\geq 0\,:\,X_{\tau_\ell+t}=v_{\ell-s}\}.\]
Notice that if $\tau_{\ell}\leq n< \tau_\ell+\tau'_{\ell,s}$, then $[T_n,X_n]\in O(s)$ almost surely. Indeed, $\tau_{\ell}\leq n< \tau_\ell+\tau'_{\ell,s}$ has the following consequences:
\begin{enumerate}
\item $X_n$ and $v_{\ell}$ lie on ``opposite sides"~of $v_{\ell-s}$, so ${\rm dist}_{T_n}(X_n,v_\ell)\geq {\rm dist}_{T_n}(v_{\ell-s},v_\ell)=s$.
\item The subtree consisting of $v_{\ell-s}$ and all nodes on the same side as $X_n$ is finite, as it consists of $v_i$ for $0\leq i\leq \ell-s$ plus the new nodes added between times $\tau_\ell$ and $n$.
\item The subtree consisting of $v_{\ell-s}$ and all nodes on the same side as $v_\ell$ is  a.s. infinite, as it contains the tree $T_0$ at time $0$.\end{enumerate} Combining these properties, we see that all long enough paths in $T_n$ that start from $X_n$ must pass through $v_\ell$, which lies at distance $\geq s$ from $X_n$. This implies $[T_n,X_n]\in O(s)$, since these same paths must pass through the unique vertex at distance $s$ from $X_n$ that lies on the path from $X_n$ to $v_\ell$. We have thus shown that if $\tau_{\ell}\leq n< \tau_\ell+\tau'_{\ell,s}$, then $[T_n,X_n]\in O(s)$, as claimed. 

We deduce:
\[\Pd_{P_p,p}([T_n,X_n]\in O(s))\geq \Pd_{P_p,p}(\tau_{\ell}\leq n< \tau_\ell+ \tau'_{\ell,s}).\]
Since
\[\Pd_{P_p,p}(\tau_{\ell}\leq n< \tau_\ell+ \tau'_{\ell,s}) \geq 1 - \Pd_{P_p,p}(\tau_\ell> n) - \Pd_{P_p,p}(\tau'_{\ell,s}<+\infty),\]
we may observe as in the proof of Lemma \ref{lem:forget} that
\[\Pd_{P_p,p}(\tau'_{\ell,s}<+\infty)\leq Ce^{-\beta\,(\ell-s)}\leq \frac{\varepsilon}{2}\]
if $\ell$ is large enough. Since $P_p$ is $p$-escapable, we may ensure that $ \Pd_{P_p,p}(\tau_\ell> n)\leq \frac{\epsilon}{2}$ by taking $n$ large enough (in terms of $\ell$). So for some choice of $\ell,n$ we obtain the desired inequality:
\[\Pd_{P_p,p}(\tau_{\ell}\leq n< \tau_\ell+ \tau'_{\ell,s}) \geq 1- \varepsilon.\]\end{proof}

\section{Final comments}\label{sec:final}

Our analysis leaves open many problems about the BGRW process and its variants. We briefly discuss five of these. 

\begin{problem}Is the speed $c(p)$ given by Theorem \ref{thm:positivespeed} a monotone function of $p$?\end{problem}

It is possible to show using our techniques that $p\mapsto c(p)$ is continuous. Monotonicity is much less clear, though our (limited) simulations suggest it should hold.

\begin{problem}Give bounds on $c(p)$ as a function of $p$.\end{problem}

For small $p$, one can extract from our proof a lower bound of the form $c(p)\geq \exp(-Cp^{-C})$ for some universal $C>0$. The best upper bound we know is $c(p)\leq p$, which comes from noting that, if $T_0$ is finite, then $T_n$  has $(1+o(1))\,pn$ vertices a.s.. 

\begin{problem}Prove more properties of the random measure $P_p$.\end{problem}

All we know about this measure is what is contained in Theorem \ref{thm:oneended2}. 

\begin{problem}Consider a model where $p=p_n$ decreases with $n$, and find a threshold function for transience.\end{problem}

That is, how fast can $p_n$ decay while maintaining transience? Clearly, our arguments break down in this regime.

\begin{problem}Consider models with cycles.\end{problem}

A ``cheap" way to create cycles would be to connect the current vertex to some uniformly random vertex at distance $K$, with $K$ a random variable. If $K$ has light tails, it should still be possible to use the local topology to study the structure of our model.

\noindent {\sc Acknowledgement.} We thank Remco van der Hofstad for suggesting the one-endedness statement in Theorems \ref{thm:oneended} and \ref{thm:oneended2}.

\appendix 

\section{A simple lemma on dependent indicators}\label{sec:appendixstochasticdomination}

We prove here Lemma \ref{lem:stochdom}.

\begin{proof}[Proof of  Lemma \ref{lem:stochdom}]For $r=1,\dots,\lfloor m/k \rfloor$ define blocks 
events 
\[F_r:=\{\forall j=(r-1)\,k+1,(r-1)\,k+2\,\dots, rk\,:\,I_j=1\}.\]
Note that $k$ consecutive indices $j$ appear in each $F_r$, so:
\[\Pd\left(\mbox{ at least $k$ consecutive $1$'s in the sequence $(I_j)_{j=1}^m$}\right)\geq \Pd\left(\bigcup_{r=1}^{\lfloor m/k\rfloor}F_r\right).\]
Moreover, the values of indices $j$ involved in events $F_r,F_{r'}$ are disjoint for $r\neq r'$. With this in mind one may easily show (via our assumptions) that $\Pd(F_1)\geq \mu^k$ and $\Pd(F_r\mid F_1^c\cap \dots F_{r-1}^c)\geq \mu^k$. 
\[1 - \Pd\left(\bigcup_{r=1}^{\lfloor m/k\rfloor}F_r\right) = \Pd\left(\bigcap_{r=1}^{\lfloor m/k\rfloor}F_r^c\right)\leq (1-\mu^k)^{\lfloor m/k\rfloor}.\]
\end{proof}
\section{Some facts on weak convergence and the local topology}\label{sec:weaklocal}

\subsection{General facts} 

In this section, $(M,d)$ is a Polish metric space and $\Prob(M,d)$ is the set of all probability measures over the Borel $\sigma$-field over $M$. We need criteria to ascertain the weak convergence of a sequence $\{P_n\}_{n\in\N}\subset \Prob(M,d)$.

\subsubsection{Criterion for tightness}
\begin{proposition}\label{prop:doubleindex}Suppose we are given a doubly-indexed family
\[F^{(r)}_{k}\,:\,r\in\N,\,k\in\N\]
of closed subsets of $M$ with the following two properties. 
\begin{enumerate}
\item For each $r\in\N$, the sets $F_k^{(r)}$ are increasing and satisfy
\[F_k^{(r)}\nearrow M\mbox{ as }k\to +\infty.\]
\item For all choices of $k_1,k_2,k_3,\dots\in\N$, the set
\[K(\{k_r\}_{r=1}^{+\infty}):=\bigcap_{r=1}^{+\infty}\,F_{k_r}^{(r)}\]
is compact.
\end{enumerate}
Assume $\{P_n\}_{n\in\N}\subset \Prob(M,d)$ is such that 
\[\forall r\in\N\,:\, \lim_{k\to +\infty}\liminf_{n\in\N}P_n(F^{(r)}_k)=1.\]
Then $\{P_n\}_{n\in\N}$ is tight.\end{proposition}
\begin{proof}Given $\epsilon>0$, we will argue that there exists a choice of sequence $\{k_r\}_{r=1}^{+\infty}$ for which
\[\mbox{\bf Goal: } \forall n\in\N\,:\, P_n(K(\{k_r\}_{r=1}^{+\infty}))\geq 1-\e.\]
Notice that, no matter what sequence we choose, we get:
\[\forall n\in\N\,:\, 1 - P_n(K(\{k_r\}_{r=1}^{+\infty})) = P_n\left(\bigcup_{r=1}^{+\infty}\,(F^{(r)}_{k_r})^c\right)\leq \sum_{r=1}^{+\infty}P_n((F^{(r)}_{k_r})^c).\]
So it suffices to show that for each $r\in\N$ we can choose $k_r\in\N$ with $\sup_nP_n((K^{(r)}_{k_r})^c)\leq \epsilon/2^r$.

To do this, fix some $r$. By our assumption that $\lim_{k\to +\infty}\liminf_{n\in\N}P_n(F^{(r)}_k)=1$, we can find $k^*_r\in\N$, $n_0\in\N$ such that\[\forall n\geq n_0\,:\,P_n(F^{(r)}_{k^*_r})\geq 1-\epsilon/2^r.\] On the other hand, for $n<n_0$ we still have:
\[P_n(F^{(r)}_{k})\nearrow 1\mbox{ as }k\to +\infty\mbox{ due to assumption 1}.\]
Since there are only finitely many $n<n_0$, we can find $k^{**}_r$ such that
\[\forall n<n_0\,:\, P_n(F^{(r)}_{k^{**}_r})\geq 1- \frac{\epsilon}{2^r}.\]
We may now choose $k_r:=\max\{k^*_r,k^{**}_r\}$. Since the sets $F^{(r)}_k$ are nested, we see at once that
\[\forall n\in\N\,:\, P_n(F^{(r)}_{k_r})\geq 1- \frac{\epsilon}{2^r},\]
which is the property we needed.\end{proof}

\subsubsection{Criterion for convergence}

\begin{proposition}\label{prop:convergenceappendix}Assume $\{P_n\}_{n\in\N}$ is a tight sequence of probability measures, with a sequence of compact sets $\{K_k\}_{k\in\N}$ attesting tightness in the sense that. 
\[\lim_k\{\liminf_{n}P_n(K_k)\}=1.\]
Then there exists a countable family $\sA$ of bounded Lipschitz functions from $M$ to $\R$, which depends only on $\{K_k\}_{k\in\N}$, such that, if $\lim_nP_nf$ exists for all $f\in\mathcal{A}$, then $P_n\Rightarrow P$ for some probability measure $P$. In that case, $P$ is uniquely characterized by the fact that $Pf=\lim_n P_nf$ for all $f\in\sA$.\end{proposition}

\begin{proof}It is known that a tight sequence $\{P_n\}_{n\in\N}$ converges weakly if and only $\lim_nP_nf$ exists for each function in the set
\[\mathcal{A}_0:=\{f:M\to\R\,:\,\|f\|_{\infty}\leq 1\mbox{ and }\|f\|_{\rm Lip}\leq 1\}.\] In that case, $Pf$ is uniquely defined by the fact that $Pf = \lim_nP_nf$ for all $f\in\mathcal{A}_0$. 

Now consider a positive rational $\epsilon>0$. The fact that $\{P_n\}_{n\in\N}$ is tight implies that there exists a compact set $K_k\subset M$, with $k=k(\epsilon)$, so that with $P_n(K_\epsilon)\geq 1-\epsilon$ for all large enough $n$. For each such set, the functions $\mathcal{A}_0$ restricted to $K_{k(\epsilon)}$ form a bounded equicontinuous family. Thus the Ascoli-Arz\`{e}la Theorem implies that there exists a countable subset $A_{\epsilon}\subset \mathcal{A}_0$ with
\begin{equation}\label{eq:deltaepsilonf}\forall \delta>0\,\forall f\in \mathcal{A}_0\,\exists f_{\epsilon,\delta}\in A_{\epsilon} \,:\, \sup_{x\in {k(\epsilon)}}|f(x)-f_{\epsilon,\delta}(x)|\leq \delta.\end{equation}
We {\em claim}
\[\mathcal{A}:=\bigcup_{\epsilon \in\Q_+}A_{\epsilon}\]
is the set we are looking for. For assume that $\lim_nP_nf$ exists for each $f\in \mathcal{A}$. For each $f\in\mathcal{A}_0$, one can choose $f_{\epsilon,\delta}\in \mathcal{A}$ as in (\ref{eq:deltaepsilonf}) and obtain:
\[|P_n(f - f_{\epsilon,\delta})|\leq |P_n\Ind{K_{k(\epsilon)}}(f - f_{\epsilon,\delta})| + P_n(K_{k(\epsilon)}^c)\leq  \delta+\epsilon\mbox{ for large enough $n$}.\]
With this, one can show that $\{P_nf\}_{n\in\N}$ is Cauchy. Indeed, 
\[\limsup_{m,n\to +\infty}|(P_n - P_m)\,f|\leq 2(\epsilon + \delta) +  \limsup_{m,n\to +\infty}|(P_n - P_m)\,f_{\epsilon,\delta}| = 2(\epsilon+\delta)\]
since we are assuming $\{P_nf_{\epsilon,\delta}\}_n$ is Cauchy. Letting $\epsilon,\delta\to 0$, we obtain that $\{P_nf\}_n$ is Cauchy, as desired. A similar reasoning shows that:
\[\lim_{n}P_nf = \lim_{\epsilon,\delta\to 0}\lim_nP_nf_{\epsilon,\delta}\]
exists for all $f\in\mathcal{A}_0$. So $P_n\Rightarrow P$. Passing to limits, we obtain from the above estimates that
\[|P(f - f_{\epsilon,\delta})|\leq \epsilon+\delta.\]
In particular, $P$ is uniquely defined by the values of $Pf=\lim_nP_nf$ for $f\in\mathcal{A}_0$, which are uniquely specified by the numbers $Pf = \lim_nP_nf$ for $f\in \mathcal{A}$.\end{proof} 

\subsection{Application to local topology over rooted graphs}\label{sec:criterionlocalweak}

The general theory of the previous section will now be applied to the space $\sT_*$ of rooted, locally finite trees. 

\begin{proof}[Proof of Proposition \ref{prop:convergencelocal}] We apply the criteria for tightness and convergence in Propositions \ref{prop:doubleindex} and \ref{prop:convergenceappendix}. For each $r,k\in\N$ the set:
\[F^{(r)}_k:=\{[T,o]\in\sT_*\,:\,D_r([T,o])\leq k\}\] is closed and $\cup_kF^{(r)}_k=\sT_*$ for each fixed $r$. Lemma 3.5 in \cite{Bordenave2016} implies that all sets of the form $K(\{k_r\}_{r\in\N})=\cap_{r=1}^{+\infty}F^{(r)}_{k_r}$ are compact. Moreover, Assumption $1$ implies that \[\forall r\in\N\,:\,\lim_k\liminf_nP_n(F^{(r)}_k)=1.\] So $\{P_n\}_{n\in\N}$ is tight by virtue of Proposition \ref{prop:doubleindex}.

We now use Proposition \ref{prop:convergenceappendix} to show convergence. That Proposition gives us a countable family of functions such that if $P_nf$ converges for all $f\in\mathcal{A}$, then $P_n$ converges. In the present setting, we need a countable family of $r$-local functions. To do this, 
define for each $n\in\N$ the projection map $\Pi_m:\sT_*\to \sT_*$ that takes $[G,o]$ to the $m$-neighborhood $[G,o]_m$.  This map is a contraction of $\sT_*$ and is therefore continuous, and moreover 
\[\forall [G,o]\in \sG_*\,:\, d([G,o],\Pi_m([G,o])) = d([G,o],[G,o]_m)\leq \frac{1}{m+1}.\]
In particular, this implies that for any bounded Lipschitz $f:\sT_*\to\R$,
\[\sup\{|f([G,o]) - f\circ \Pi_m([G,o])|\,:\,[G,o]\in K\} \leq \frac{\|f\|_{\rm Lip}}{m+1}.\]
It transpires that one can replace the family $\mathcal{A}$ with
\[\mathcal{S}:=\{f\circ \Pi_m\,:\, f\in\mathcal{A},m\in\N\},\]
which is still countable, and obtain the desired result. \end{proof}

\bibliographystyle{plain}
\bibliography{ref}

\end{document}